\documentclass[preprint,9pt]{elsarticle}
\usepackage{amssymb}
\usepackage{amsmath}
\usepackage{amsthm}
\usepackage{hyperref}
\usepackage{amsfonts}
\usepackage{color}
\usepackage{epsfig}
\usepackage{graphics,subfigure}                 
\usepackage{graphicx}
\usepackage{float}
\usepackage{epsf}
\usepackage{mathrsfs}
\usepackage{multirow}

\newtheoremstyle{thmm}{1.5ex plus 1ex minus .2ex}{1.5ex plus 1ex minus
.2ex}{\rmfamily}{}{\bfseries}{}{1em}{} \theoremstyle{thmm}
\newtheorem{theorem}{Theorem}[section]
\newtheorem{lemma}{Lemma}[section]

\newtheorem{example}{Example}[section]

\biboptions{numbers,sort&compress}
\allowdisplaybreaks

\def\refe#1{(\ref{#1})}

\def\curl{\mathrm{curl}}
\def\div{\mathrm{div}}
\def\u{\boldsymbol{u}}
\def\b{\boldsymbol{b}}
\def\v{\boldsymbol{v}}
\def\c{\boldsymbol{c}}
\def\w{\boldsymbol{w}}
\def\x{\boldsymbol{x}}
\def\q{\boldsymbol{q}}
\def\m{\boldsymbol{m}}

\def\n{\boldsymbol{n}}
\def\f{\boldsymbol{f}}
\def\g{\boldsymbol{g}}
\def\ps{\boldsymbol{\psi}}
\def\Ph{\boldsymbol{\Phi}}
\def\lmd{\lambda}

\def\gm{\gamma}
\def\ds{\displaystyle}

\def\H{\boldsymbol{H}}
\def\M{\boldsymbol{M}}
\def\Z{\boldsymbol{Z}}
\def\I{\mathrm{I}}

\def\V{\boldsymbol{V}}
\def\W{\boldsymbol{W}}
\def\T{\mathcal{T}}
\def\S{\mathrm{S}}
\def\O{\Omega}
\def\Ha{\mathrm{Ha}}

\begin{document}
\begin{frontmatter}
\title{\bf Conservative nonconforming virtual element method for stationary incompressible magnetohydrodynamics 
\thanks{.}}

\author[]{Xiaojing Dong}
\ead{dongxiaojing99@xtu.edu.cn}

\author[]{Yunqing Huang}
\ead{huangyq@xtu.edu.cn}

\author[]{Tianwen Wang\corref{mycorrespondingauthor}}
\cortext[mycorrespondingauthor]{Corresponding author}
\ead{wangtianwen@smail.xtu.edu.cn}

\address{Hunan Key Laboratory for Computation and Simulation in Science and Engineering, Key Laboratory of Intelligent Computing $\&$ Information Processing of Ministry of Education, School of Mathematics and Computational Science, Xiangtan University, Xiangtan, Hunan, 411105, P.R. China}

\begin{abstract}

In this paper, we propose a conservative nonconforming virtual element method for the full stationary incompressible magnetohydrodynamics model. We leverage the virtual element  satisfactory divergence-free property to ensure mass conservation for the velocity field.  The condition of the well-posedness of the proposed method, as well as the stability are derived.  We  establish optimal error estimates in the discrete energy norm for both the velocity and magnetic field. Furthermore, by employing a new technique, we obtain the optimal error estimates in $L^2$-norm without any additional conditions. Finally, numerical experiments are presented to validate the theoretical analysis. In the implementation process, we adopt the effective Oseen iteration to handle the nonlinear system.


\noindent{ \bf Keywords:} Virtual element method; Stationary incompressible magnetohydrodynamics; error estimates; Polygonal meshes
\end{abstract}
\end{frontmatter}

\section{Introduction}
Incompressible magnetohydrodynamics (MHD) equations, coupled by multiple physical fields, model the interaction between conducting fluid and an external magnetic field. MHD equations are significant in the regions of physics and engineering, such as cooling of liquid metal in nuclear reactors, confinement for controlled thermonuclear fusion, magnetic fluid motors and the propagation of radio waves in the ionosphere \cite{MHD-2006,MHD-2017}. 

In the last decades, there has been a lot of research on numerical solutions of incompressible MHD equations (see \cite{G-MHD-1991,MHD-2006,Mixed-MHD,FullDivFree-MHD,Huang-DG,Dong-Three,Dong-Two-level,Hu-Stationary,Hu-DivFree-B,Tang-MHD} and the references therein). Regarding handling of nonlinear terms, the classic finite element iterative methods were designed in \cite{Dong-Three,Dong-Two-level}. In the most mixed finite element methods, the divergence-free restriction is imposed in the weak sense and the numerical solutions violate the physical properties. The divergence-free methods have been received widely attention for incompressible MHD equations. A mixed finite element method \cite{Mixed-MHD} (FEM) and three interior penalty discontinuous Galerkin (DG) methods \cite{Huang-DG} for incompressible MHD were developed, with preserving exactly divergence-free velocity. Literatures \cite{Hu-Stationary,Hu-DivFree-B} studied structure-preserving FEMs of MHD. Moreover, a fully divergence-free FEM for MHD was proposed in \cite{FullDivFree-MHD}, where the magnetic field is approximated by a magnetic vector potential.

Virtual element method (VEM), introduced with one of the simplest examples in \cite{Basis-VEM,Guide-VEM}, is a recent technique for numerical approximation of partial differential equations and is regarded as the extension of classic FEM onto polygonal/polyhedral meshes. VEM does not require that shape functions on element are explicitly constructed, but the discrete algebraic systems can be obtained by using the degrees of freedom. An arbitrary order nonconforming virtual element was presented and applied to the Poisson equation \cite{NCVEM-Poisson} and the Stokes equations \cite{NCVEM-Stokes,Liu-Stokes}. The literature \cite{DivFree-Stokes} designed an $\H^1$-conforming Stokes-type virtual element $(k\geq 2)$ such that the associated discrete kernel is divergence-free. In \cite{Stream-Stokes,Stream-NS}, the stream virtual element formulations were studied by using a suitable stream function space to approximate velocity. Piecewise divergence-free nonconforming virtual elements for Stokes problem in any dimensions were investigated in \cite{Wei-Stokes}. Based on the de Rham complex, the reference \cite{Zhao-Stokes} constructed a nonconforming virtual element for Stokes problem, which is piecewise divergence-free in discrete kernel. When utilizing the virtual element proposed in \cite{Zhao-Stokes} to solve the Navier-Stokes equations, to ensure that the constructed discrete VEM scheme achieves optimal convergence with suitable precision, the computable $L^2$-projection onto the polynomials of degree $\leq k$ is required. In order to calculate $L^2$-projection, the enhanced version based on the nonconforming virtual element \cite{Zhao-Stokes} was presented in \cite{Zhao-NS} for the Navier-Stokes equations. In order to derive optimal convergence with suitable precision, the enhanced version of the $H^1$-conforming virtual element proposed in \cite{Basis-VEM} was investigated in \cite{Equivalent-VEM} for  the reaction-diffusion problems. VEM has additionally been used to solve the resistive MHD model \cite{VEM-MHD-2D,VEM-MHD-3D}. 

Motivated by \cite{VEM-MHD-2D,VEM-MHD-3D}, we consider a nonconforming VEM for the 2D full stationary incompressible MHD equations. Compared to the work in  \cite{VEM-MHD-2D,VEM-MHD-3D}, our research has two main differences: first, the model discussed here is more complex than the resistive MHD; second, the study in  \cite{VEM-MHD-2D,VEM-MHD-3D} focuses on low-order virtual elements, while we develop an arbitrary order ($k\geq 1$) scheme on shape-regular polygonal meshes. For keeping mass conservation, the velocity and the pressure are approximated by the enhanced nonconforming virtual element space and discontinuous piecewise polynomials. The discrete velocity is piecewise divergence-free. Each component of the magnetic field is approximated by the enhanced $H^1$-conforming virtual element space. The treatment approaches for the bilinear term related to the magnetic field and the complex trilinear terms are provided. We prove the stability of the discrete formulation and establish conditions for uniqueness of the solution. Additionally, we derive optimal error estimates in the discrete energy norm and $L^2$-norm for both velocity and magnetic field, as well as in the $L^2$-norm for the pressure.

The rest of the paper is organized as follows. In Section 2, we briefly introduce some basis notations, the full stationary incompressible MHD equations and the weak formulation. In Section 3, we present the approximate spaces of the velocity, the magnetic field and the pressure, the virtual element discretizations of the MHD equations and the stability of the discrete formulation. In Section 4, the optimal error estimates are derived. Section 5 shows some numerical experiments to validate the theoretical analysis.

\section{Preliminaries}
Let $\O \subset \mathbb{R}^2$ be a bounded convex polygonal domain. We denote a two-dimensional variable by $\x$, where $\x= (x_1,x_2)$ and $\x^{\perp} = (-x_2,x_1)$. We use $\n $ to represent the unit outward normal vector of $\partial \O$. For any polygonal subdomain $E$ of $\O$, the unit outward normal vector of $\partial E$ (edge $e$ of $\partial E$) is denoted $\n_E \,( \n_e)$. We adopt $\mathbb{P}_k(\mathcal{O})$ to express the set of polynomials on $\mathcal{O}$ of degree at most $k(\ge 0)$, here $\mathcal{O}$ can be $E$, $\partial E$ or edge $e$, specially $\mathbb{P}_{-2}(\mathcal{O}) = \mathbb{P}_{-1}(\mathcal{O}) = \{\emptyset \}$. There exists a direct sum decomposition
\begin{align*}
\begin{aligned}
\left[ \mathbb{P}_k(E) \right]^2 &= \nabla \mathbb{P}_{k+1}(E) \oplus \x^{\perp} \mathbb{P}_{k-1}(E).
\end{aligned} 
\end{align*}
For scalar function $v$ and vector functions $\v = (v_1,v_1)$, $\w = (w_1,w_2)$, we introduce the cross product and $\curl$ operators
\begin{align*}
\begin{aligned}
v \times \w &= (- w_2 v, w_1 v ) , \qquad   \v \times \w = w_2 v_1 - w_1 v_2,\\
\curl\, v &= \left( \frac{\partial v}{\partial x_2}, -\frac{\partial v}{\partial x_1} \right),  \qquad \curl \,\v =  \frac{\partial v_2}{\partial x_1}  - \frac{\partial v_1}{\partial x_2} .
\end{aligned}
\end{align*}

We define some useful spaces as follows:
\begin{align*}
&H^1(\O) = \left\{ v \in L^2(\O) : \; \nabla v \in L^2(\O) \right\},\qquad H_0^1(\O) = \left\{ v \in H^1(\O) :  v|_{\partial \O} = 0 \right\},\\
&\H^1(\O) = \left[ H^1(\O) \right]^2,\qquad
\H_0^1(\O) = \left[ H_0^1(\O) \right]^2,\qquad
\boldsymbol{L}^2(\O) = \left[ L^2(\O) \right]^2,\\
&\H_n^1 = \left\{ \v \in \H^1(\O) : \v\cdot \n|_{\partial \O } = 0 \right\}, \qquad  L_0^2(\O) = \left\{ q \in L^2(\O) : \; \int_{\O} q \,\mathrm{d} \O = 0 \right\}.
\end{align*}
Let $W^{k,p}(E)$ be the standard Sobolev space equipped with norm $\| \cdot\|_{k,p,E}$ and seminorm $|\cdot|_{k,p,E}$. In particular, $W^{0,p}(E)$ is Lebesgue space with norm $\| \cdot\|_{L^p(E)}$, $W^{k,2}(E)$ is Hilbertian space with norm $\| \cdot\|_{k,E}$ and seminorm $|\cdot|_{k,E}$. $(\cdot,\cdot)_{E}$ is usual $L^2$-inner product. If $E = \O$, the subscript $E$ of norm, seminorm and inner product will be omitted. 

We consider the two-dimensional full stationary incompressible MHD \cite{G-MHD-1991}:
\begin{align}
\begin{cases}
-R_{\nu}^{-1} \Delta \u + (\nabla \u)\u  + \nabla p - S_c \,\curl\,\b \times \b = \f ,\quad \mathrm{in} \; \O,  \\
R_m^{-1} S_c\,\curl\,(\curl\,\b) - S_c\,\curl\,(\u \times \b) = \g, \quad \mathrm{in} \; \O, \\
\div\,\u = 0, \quad \mathrm{in} \; \O,\\
\div\,\b = 0, \quad \mathrm{in} \; \O,\\
\u = \mathbf{0} , \quad \mathrm{on} \; \partial \O,\\
\b\cdot \n = 0, \;\n\times \curl\,\b = \mathbf{0}, \quad \mathrm{on} \; \partial \O,
\end{cases}    \label{MHDeq}
\end{align}
where $\u,\,\b,\,p$ are the velocity, the magnetic field and the pressure, respectively. Hydrodynamic Reynolds number $R_{\nu}$, magnetic Reynolds number $R_m$ and coupling coefficient $S_c$ are the physical parameters of the equations. Functions $\f$ and $ \g \in \boldsymbol{L}^2(\O)$ are source terms.

Like in \cite{G-MHD-1991}, the weak formulation of problem \refe{MHDeq} is as follows:  $\forall (\v,\c,q) \in \H_0^1(\O)\times \H_n^1(\O)\times L_0^2(\O) $, find $(\u,\b,p) \in \H_0^1(\O) \times \H_n^1(\O)\times L_0^2(\O)$ such that 
\begin{align}
\begin{cases}
R_{\nu}^{-1} (\nabla \u,\nabla \v) + \frac{1}{2}((\nabla \u)\u,\v) - \frac{1}{2}((\nabla \v)\u,\u) - (\div\,\v,p) - S_c\, (\curl \,\b\times \b,\v) = (\f,\v), \\
 R_m^{-1} S_c\,(\curl \,\b,\curl\, \c) +  R_m^{-1} S_c\,(\div\,\b,\div\, \c)+ S_c\,(\curl\,\c\times \b,\u) = (\g,\c) ,\\
(\div\,\u,q) = 0.
\end{cases} \label{MHDweakeq}
\end{align}
According to the Proposition 3.1 in \cite{G-MHD-1991}, we know that $\div\, \b = 0$ can be derived from \refe{MHDweakeq}. Problem \refe{MHDweakeq} can be rewritten as follows: find $(\u,\b,p) \in \H_0^1(\O) \times \H_n^1(\O)\times L_0^2(\O)$ such that
\begin{align}
\begin{cases}
A(\u,\b;\v,\c) + B(\u,\b;\u,\b;\v,\c)  - d(\v,\c;p)   = \left< F;\v,\c\right>,\quad \forall \v \in \H_0^1(\O), \c\in  \H_n^1(\O), \\
d(\u,\b;q) = 0, \quad \forall q\in L_0^2(\O),
\end{cases}   \label{MHDrewriteeq}
\end{align} 
where, $\forall \w\in \H_0^1(\O), \ps\in\H_n^1(\O)$,
\begin{align*}
& A(\u,\b;\v,\c) = a_0(\u,\v) + a_1(\b,\c),\\ 
& B(\u,\b;\w,\ps;\v,\c) = a_2(\boldsymbol{\u},\w,\v) - a_3(\ps,\b,\v) + a_3(\c,\b,\w),\quad   \\
& a_0(\u,\v) = R_{\nu}^{-1}(\nabla \u,\nabla \v) ,\\
& a_1(\b,\c) = R_m^{-1}S_c\, (\curl \,\b,\curl\, \c) + R_m^{-1}S_c\, (\div\,\b,\div\, \c),\\
& a_2(\u,\w,\v) = \frac{1}{2}((\nabla \w)\u,\v) - \frac{1}{2}((\nabla \v)\u,\w), \\
& a_3(\ps,\b,\v) = S_c\,(\curl\,\ps\times\b,\v),\\
& d(\v,\c;q) = (\div\,\v,q), \quad \left<F;\v,\c\right> = \left(\f,\v\right)  + (\g,\c).
\end{align*}   
For the sake of convenience, some norms are defined by   
\begin{align*}
& \|(\v,\c)\|_{i} = \left( \| \v\|_i^2 + \| \c\|_i^2 \right)^{\frac{1}{2}},\qquad \forall \v,\c \in \H^i(\O), \; i = 0,1,2,
\end{align*}
and $\|F\|_{\ast} = \|(\f,\g)\|_0$. We use the equivalent norm $\|\nabla \v\|_0$ to replace $\|\v\|_1$ for all $ \v \in \H_0^1(\O)$. From \cite{NS-1986}, it know that
\begin{align}
& \|\u \|_{L^4} \le \lmd_0 \| \u\|_1 , \qquad  \| \curl \,\b \|_0^2 + \|\div\,\b \|_0^2  \geq  \lmd_1\|\b\|_1^2, \qquad \forall \u \in \H^1(\O) ,\; \b \in \H_n^1(\O). \label{Soboleveq}
\end{align}  
Here and after, $\lmd_i,i=0,1,...,5$ are positive constants independent of mesh size. 

For the above linear forms, the following estimates hold (see \citep{G-MHD-1991}),  $\forall \u,\w,\v \in \H_0^1(\O)$, $\b,\ps,\c \in \H_n^1(\O)$ :
\begin{align}
|A(\u,\b;\v,\c)| & \le \max \left\{ R_{\nu}^{-1},4R_m^{-1}S_c \right\} \|(\u,\b) \|_1 \,\|(\v,\c)\|_1, \label{cMHDeq1}\\
A(\u,\b;\u,\b) & \ge \min\{R_{\nu}^{-1}, \lmd_1R_m^{-1}S_c \} \|(\u,\b)\|_1^2,\label{cMHDeq2}\\
|B(\u,\b;\w,\ps;\v,\c)| & \le \sqrt[]{2}\lmd_0^2 \max\{1,\sqrt[]{2} S_c \} \|(\u,\b) \|_1 \, \|(\w,\ps) \|_1 \,\|(\v,\c) \|_1. \label{cMHD3}
\end{align} 
The linear form $d(\cdot,\cdot;\cdot)$ is continuous and satisfies the inf-sup condition \cite{NS-1986,G-MHD-1991}
\begin{align}
\sup_{(\v,\c) \in \H_0^1(\O) \times \H_n^1(\O)}  \frac{|d(\v,\c;q)|}{\|(\v,\c) \|_1} \ge \gm_0 \| q \|_0,\qquad \forall \, q \in L_0^2(\O) . 
\end{align}   

\begin{theorem}\label{c-exist} (\cite{Dong-Three})
 If $\frac{\sqrt[]{2} \lmd_0^2 \max\{1,\sqrt[]{2} S_c\}\|F\|_{\ast}}{(\min\{R_{\nu}^{-1},\lmd_1 R_m^{-1}S_c\})^2} < 1$, the problem \refe{MHDrewriteeq} has a unique solution $(\u,\b,p) \in \H_0^1(\O) \times \H_n^1(\O) \times L_0^2(\O) $, satisfying 
\begin{align}
\min\{ R_{\nu}^{-1},\lmd_1 R_m^{-1}S_c\} \|(\u,\b)\|_1 \le \|F\|_{\ast}.
\end{align} 
\end{theorem}
\setcounter{equation}{0}
\section{Virtual element method}
\subsection{Meshes and projections}
Let $\{\T_h\}_h$ be a family of partitions of $\O$ into convex polygons and $h_E$ be diameter of element $E$, $h$ be the largest of all diameters, $\mathcal{E}_h$ represent the set of edges of $\T_h$. For each element $E$, we have mesh regular assumptions that there exists a positive constant $\rho $ satisfying:
\begin{flushleft}
  (\textbf{A1}) $E$ is star-shaped with respect to a ball $B_E$ of radius $\ge \rho h_E$.  \\
  (\textbf{A2}) The distance between any two vertexes of $ E$ is $\geq \rho h_E$.
\end{flushleft}
For an internal edge $e$, which is a shared edge of elements $E^{+}$ and $E^{-}$ of $\T_h$, we define the jump of $\v$ by
$[\v]|_e = (\v|_{\partial E^+})|_e - (\v|_{\partial E^-})|_e$, specially, for boundary edge $e$, $[\v]|_e = \v|_e$.

We introduce some useful projections as follows:
\begin{flushleft}
 \begin{itemize}
   \item $H^1$-seminorm projection $\Pi_k^{\nabla,E}$ : $\H^1(E) \rightarrow [\mathbb{P}_k(E)]^2$, defined by
   \begin{align}
   \begin{cases}
    \int_E \nabla \Pi_k^{\nabla,E} \v : \nabla \m_k \,\mathrm{d} E = \int_E \nabla \v : \nabla \m_k \,\mathrm{d}E, \qquad \forall \m_k \in \left[ \mathbb{P}_k(E) \right]^2, \;\v \in \H^1(E),\\
    \int_{\partial E} \Pi_k^{\nabla,E} \v \,\mathrm{d} \partial E = \int_{\partial E} \v \;\mathrm{d}\partial E ,\qquad \forall \v\in \H^1(E).
   \end{cases}    \label{H1-projection}
   \end{align}
   \item $L^2$-projection $\Pi_k^{0,E}$ : $\boldsymbol{L}^2(E) \rightarrow  \left[\mathbb{P}_k(E)\right]^2 $, defined by
   \begin{align}
   \int_E \Pi_k^{0,E} \v \cdot \m_k \,\mathrm{d}E = \int_E\v \cdot \m_k \,\mathrm{d}E ,  \qquad \forall \m_k \in \left[ \mathbb{P}_k(E) \right]^2,\, \v \in \boldsymbol{L}^2(E). \label{L2-projection}
   \end{align}
 \end{itemize}
\end{flushleft}

\subsection{The nonconforming element space and piecewise polynomials space}
We use the enhanced nonconforming virtual element \cite{Zhao-NS} to approximate the velocity. In order to describe it, we first introduce two auxiliary spaces. One space is
\begin{align}
\V_{k}^{f}(E) = \left\{ \v \in \H^1(E) : \; \div \,\v \in \mathbb{P}_{k-1}(E),\; \curl \,\v \in \mathbb{P}_{k-1}(E), \; \v\cdot \n_e|_e \in \mathbb{P}_k(e) , \; \forall e \in \partial E   \right\}. 
\end{align}
From \cite{Zhao-Stokes}, we know that the following direct sum decomposition holds
\begin{align}
\V_k^f(E) = \boldsymbol{W}_0(E) \oplus  \boldsymbol{W}_1(E),  
\end{align} 
where $\boldsymbol{W}_0(E)$ and $\boldsymbol{W}_1(E)$ are 
\begin{align*}
\begin{aligned}
\boldsymbol{W}_0(E) &= \left\{ \v \in \V_k^f(E) : \; \div\,\v = 0, \; \v \cdot \n_e|_e = 0 ,\; \forall e \in \partial E\right\}, \\
\boldsymbol{W}_1(E) &= \left\{ \v \in \V_k^f(E) : \; \curl \,\v = 0 \right\}.
\end{aligned}  
\end{align*}
The other space is 
\begin{align*}
\Psi(E) = \left\{ \psi \in H^2(E) : \;\Delta^2 \psi \in \mathbb{P}_{k-1}(E),\; \Delta \psi |_e \in \mathbb{P}_{k-1}(e),\; \psi |_e = 0 ,\;\forall e \in \partial E  \right\}. 
\end{align*}  
From Lemma 3.1 in \cite{Zhao-NS}, it is easy to see that 
\begin{align}
\boldsymbol{\widetilde{V}}(E) = \boldsymbol{W}_1(E) \oplus \curl\, \Psi (E) = \V_k^f(E) +  \curl\, \Psi (E). 
\end{align} 
Then, the local space is defined by
\begin{align}
\begin{aligned}
\V(E) = \left\{ \v \in \boldsymbol{\widetilde{V}}(E): \, \int_E \v \cdot \m_k \,\mathrm{d}E = \int_E \Pi_k^{\nabla,E} \v \cdot \m_k \,\mathrm{d}E ,\; \forall \m_k \in  \x_E^{\perp} \mathbb{P}_{k-1}(E) /\x_E^{\perp} \mathbb{P}_{k-3}(E), \right.\\
\left. \qquad \int_e \v \cdot \n_e m_k \,\mathrm{d}e =  \int_e \Pi_k^{\nabla,E} \v \cdot \n_e m_k \,\mathrm{d}e, \quad \forall m_k\in \mathbb{P}_k(e)/\mathbb{P}_{k-1}(e), e \subseteq \partial E  \right\}, \label{localNCVEM}
\end{aligned} 
\end{align}
where $\x_E = (\x -\x_b)/h_E$, $\x_b$ is the barycenter of $E$. Obviously, $ \left[ \mathbb{P}_k(E)\right]^2 \subseteq \V(E)$. The degrees of freedom of $\V(E)$ are 
\begin{align}
& \bullet \quad\frac{1}{|e|} \int_e \v \cdot \m_{k-1} \,\mathrm{d}e , \quad \forall \m_{k-1} \in \left[\mathbb{P}_{k-1}(e)\right]^2,\, e \in \partial E \label{NCVEMdof1} ,\\
& \bullet \quad \frac{1}{|E|} \int_E \v\cdot \m_{k-2} \,\mathrm{d}E,\quad \forall \m_{k-2} \in  \left[\mathbb{P}_{k-2}(E)\right]^2 .\label{NCVEMdof2}
\end{align}
It is not difficult for us to find that the projection $\Pi_k^{\nabla,E} \circ\V(E)$ can be computed by using integration by parts and \refe{NCVEMdof1}-\refe{NCVEMdof2}. According to the definition of $\V(E)$, we know that $\div \,\v \in \mathbb{P}_{k-1}(E)$ for all $\v \in \V(E)$ and the expression of $\div \, \v$ on element $E$ can be calculated by
\begin{align}
\int_E \div\, \v \cdot m_{k-1} \,\mathrm{d} E = - \int_E \v \cdot \nabla m_{k-1} \,\mathrm{d} E + \sum_{e \in \partial E}\int_{e} \v\cdot \n_e m_{k-1} \,\mathrm{d}e, \quad \forall m_{k-1} \in \mathbb{P}_{k-1}(E).  \label{NCVEMdivveq}
\end{align}
The first term on the right-hand side of \refe{NCVEMdivveq} can be computed by \refe{NCVEMdof2}. The second term is computable, since $\left(\v\cdot \n_e\right)|_{e}  \in \mathbb{P}_k(e)$ and its expression can be uniquely determined by \refe{localNCVEM}-\refe{NCVEMdof2}. For all $\m_k \in \left[\mathbb{P}_k(E)\right]^2,\v\in \V(E)$, there holds
\begin{align}
\begin{aligned}
\int_E \v\cdot \m_{k} \,\mathrm{d}E &= \int_E \v\cdot \left( \nabla m_{k+1} \oplus \q_k \right) \,\mathrm{d}E\\
&= -\int_E \div\,\v\cdot m_{k+1} \,\mathrm{d}E + \sum_{e \in \partial E}\int_{e} \v\cdot \n_e m_{k+1} \,\mathrm{d}e + \int_E \v\cdot \q_k \,\mathrm{d}E,
\end{aligned} \label{cL2projection}
\end{align}
where $m_{k+1}\in \mathbb{P}_{k+1}(E)$, $\q_k \in \x_E^{\perp}\mathbb{P}_{k-1}(E) $. Similar to the previous discussion, the right-hand side of \refe{cL2projection} can be calculated by \refe{NCVEMdof1}-\refe{NCVEMdof2}. Thus, the projection $\Pi_k^{0,E} \circ \V(E)$ is computable.

We define $L^2$-projection $\Pi_{k-1}^{0,E} $ : $ \nabla \V(E) \rightarrow  \left[ \mathbb{P}_{k-1}(E)\right]^{2\times 2} $ by 
\begin{align}
\int_E \Pi_{k-1}^{0,E} \left( \nabla \v \right): \m_{k-1} \,\mathrm{d} E  = \int_E   \nabla \v : \m_{k-1} \,\mathrm{d} E,  \quad \forall \m_{k-1} \in \left[ \mathbb{P}_{k-1}(E) \right]^{2\times 2}.  \label{gradL2projection}
\end{align}
Indeed, the right-hand side of \refe{gradL2projection} can be calculated by integration by parts and \refe{NCVEMdof1}-\refe{NCVEMdof2}. As a result, the projection $\Pi_{k-1}^{0,E}\circ (\nabla \V(E))$ is computable. With the above preparation, the global space is defined as 
\begin{align}
\V = \left\{ \v \in \boldsymbol{L}^2(\O) : \; \v|_E \in \V(E) ,\, \int_e [\v]\cdot \m_{k-1} \,\mathrm{d}e =0  ,\quad \forall \m_{k-1} \in \left[ \mathbb{P}_{k-1}(e) \right]^2,\,  e \in \mathcal{E}_h \right\}.   \label{globalNCVEM}
\end{align}   
\begin{lemma}(\cite{Zhao-NS}) \label{NCVEM-interpolation}
Let $\v \in \H^l(E)$ with $1\le l \le k+1$, $\v_I$ be its degrees of freedom interpolation. Then, it holds that 
\begin{align*}
\|\v - \v_I\|_{0,E} + h|\v - \v_I|_{1,E} \le Ch^l|\v|_{l,E}.
\end{align*}
Throughout this article, $C$ denotes a generic positive constant independent of $h$, and might be different value at each occurrence.
\end{lemma}

As the usual framework of VEM, we define a computable local bilinear form  
$a_{0h}^E (\cdot,\cdot): \V(E)\times\V(E) \rightarrow \mathbb{R}$ by
\begin{align*}
a_{0h}^E(\u_h,\v_h) &= R_{\nu}^{-1}\Big((\nabla \Pi_k^{\nabla,E} \u_h,\nabla \Pi_k^{\nabla,E} \v_h)_{E} + \S_0^E((\I -\Pi_k^{\nabla,E})\u_h,(\I -\Pi_k^{\nabla,E})\v_h) \Big),\quad \forall \u_h,\v_h \in\V(E), 
\end{align*}
where $\S_0^E$ is a symmetric and positive definite bilinear form such that
\begin{align*}
\alpha_{0\ast} a_0^E(\v_h,\v_h) \le R_{\nu}^{-1}\S_0^E(\v_h,\v_h) \le \alpha_0^{\ast} a_0^E(\v_h,\v_h) ,\quad \forall \v_h \in \V(E)\cap\mathrm{Ker}(\Pi_k^{\nabla,E}),  
\end{align*} 
for two positive constants $\alpha_{0\ast}$ and $\alpha_0^{\ast}$ independent of $h$ and $E$. The definition of $a_{0h}^E (\cdot,\cdot)$ and \refe{H1-projection} imply 
\begin{itemize}
 \item[•] $k$-consistency: for all $\m_k \in [\mathbb{P}_k(E)]^2$ and $\v_h \in \V(E)$,
  \begin{align}
  a_{0h}^E(\m_k,\v_h) = a_0^E(\m_k,\v_h), \label{NC-consistency}
  \end{align}
 \item[•] stability: there exist two positive constants $\sigma_{0\ast}$ and $\sigma_0^{\ast}$, independent of $h$ and $E$, satisfying
 \begin{align}
  \sigma_{0 \ast} a_0^E(\v_h,\v_h)  \le a_{0h}^E(\v_h,\v_h) \le \sigma_0^{\ast} a_0^E(\v_h,\v_h),\quad \forall \v_h \in \V(E).\label{NC-stability}
 \end{align}
\end{itemize}

The pressure is approximated by discontinuous piecewise polynomials, which is defined as
\begin{align*}
Q = \left\{ q \in L_0^2(\O): \,q|_E \in \mathbb{P}_{k-1}(E) ,\quad \forall \,E \in \T_h \right\}.
\end{align*}
Let $q\in H^l(\O) \cap L_0^2(\O)$, $q_I|_E$ be the $L^2$-projection of $q|_E$ onto $\mathbb{P}_{k-1}(E)$, then $q_I \in Q$, there holds
\begin{align}
\|q - q_I\|_{0,E} +  h|q - q_I|_{1,E}\leq C h^{l}|q|_{l,E}, \quad \forall 1\leq l \leq k.  \label{L2-eq}
\end{align}
\subsection{The nodal space}
In order to compute the $L^2$-projection onto polynomials of degree $\le k$, we use the enhanced $H^1$-conforming virtual element \cite{Equivalent-VEM} to approximate each component of the magnetic field. First of all, we introduce a auxiliary space
\begin{align}
U_k^n(E) = \left\{ c \in H^1(E) :\, \Delta c \in \mathbb{P}_{k}(E),\,   c|_{\partial E} \in C^0(\partial E) ,\,c|_e \in \mathbb{P}_{k}(e) ,\quad \forall \, e \in \partial E   \right\} .
\end{align}
Then, the local enhanced space is defined as 
\begin{align}
\M(E) = \left\{ \c\in [U_k^n(E)]^2:\; (\c - \Pi_k^{\nabla,E} \c,\m)_E=0, \quad \forall \m\in [\mathbb{P}_k(E)]^2/[\mathbb{P}_{k-2}(E)]^2 \right\}, \label{localCVEM}  
\end{align}
with the degrees of freedom 
\begin{align}
& \bullet \quad \mathrm{the} \; \mathrm{values} \;  \mathrm{of} \; \c \; \mathrm{at} \; \mathrm{vertices}\;  \mathrm{of}\; E,     \label{CVEMdof1} \\
& \bullet \quad  \frac{1}{|e|}\int_e \c \cdot\m_{k-2} \, \mathrm{d}e, \quad \forall \m_{k-2} \in [\mathbb{P}_{k-2}(e)]^2 , \,e \in \partial E,\label{CVEMdof2} \\
& \bullet \quad \frac{1}{|E|} \int_E  \c \cdot  \m_{k-2} \,\mathrm{d}E ,\quad \forall \m_{k-2} \in [\mathbb{P}_{k-2}(E)]^2.  \label{CVEMdof3}
\end{align}
It is obvious that the projection $\Pi_k^{\nabla,E} \circ \M(E)$ can be computed by \refe{CVEMdof1}-\refe{CVEMdof3}. From \refe{CVEMdof3} and \refe{localCVEM}, we know that the projection $\Pi_k^{0,E} \circ \M(E)$ is computable. We define $L^2$-projection  $ \; \Pi_{k-1}^{0,E} $ : $\curl\, \M(E) \rightarrow \mathbb{P}_{k-1}(E)$ by
\begin{align}
\int_E \Pi_{k-1}^{0,E} \left( \curl \, \c \right) \cdot m_{k-1} \,\mathrm{d}E = \int_E \curl\,\c \cdot m_{k-1} \;\mathrm{d}E,  \quad \forall m_{k-1} \in \mathbb{P}_{k-1}(E). \label{curlprojecteq1}
\end{align}  
According to integration by parts and \refe{CVEMdof1}-\refe{CVEMdof3}, it is easy to see that the projection $\Pi_{k-1}^{0,E} \circ (\curl \,\M(E))$ is computable. Likewise, the projection $\Pi_{k-1}^{0,E} \circ (\div\,\M(E))$ can be calculated.

The global space is defined as
\begin{align}
\M = \left\{\c \in \H_n^1(\O):\; \c|_E \in \M(E)  ,\quad \forall E \in \T_h   \right\}.\label{globalCVEM}
\end{align}
\begin{lemma}(\cite{Small-edge}) \label{CVEM-interporlation} 
Let $\c \in \H^l(E)$ with $1 \le l \le k+1$, $\c_I$ be its degrees of freedom interpolation. Then, it holds that
\begin{align*}
\|\c - \c_I\|_{0,E} + h|\c - \c_I|_{1,E} \le Ch^l|\c|_{l,E}.
\end{align*}
\end{lemma} 
 We define a computable local bilinear form  
$a_{1h}^E (\cdot,\cdot): \M(E)\times\M(E) \rightarrow \mathbb{R}$ by
\begin{align*}
a_{1h}^E(\b_h,\c_h) &= R_m^{-1}S_c\left( (\Pi_{k-1}^{0,E}\curl\, \b_h,\Pi_{k-1}^{0,E}\curl\, \c_h)_{E} + (\Pi_{k-1}^{0,E}\div\, \b_h,\Pi_{k-1}^{0,E}\div \, \c_h)_{E} \right.\\
& \left. \quad + \S_1^E((\I -\Pi_k^{\nabla,E})\b_h,(\I -\Pi_k^{\nabla,E})\c_h)\right),\quad \forall \b_h,\c_h \in\M(E), 
\end{align*}
where $\S_1^E$ is a symmetric and positive definite bilinear form satisfying
\begin{align*}
\alpha_{1\ast} \|\nabla \c_h\|_{0,E} \le \S_1^E(\c_h,\c_h) \le \alpha_1^{\ast} \|\nabla \c_h \|_{0,E},\quad \forall \c_h \in \M(E)\cap\mathrm{Ker}(\Pi_k^{\nabla,E}),  
\end{align*} 
for two positive constants $\alpha_{1\ast}$ and $\alpha_1^{\ast}$ independent of $h$ and $E$.
\begin{lemma}
The local bilinear form $a_{1h}^{E}(\cdot,\cdot)$ satisfies the following properties
\begin{itemize}
 \item[•] $k$-consistency: for all $\m_k \in [\mathbb{P}_k(E)]^2$ and $\c_h \in \M(E)$,
  \begin{align}
  a_{1h}^E(\m_k,\c_h) = a_1^E(\m_k,\c_h), \label{CVEM-cosistency}
  \end{align}
 \item[•] stability: there exist two positive constants $\sigma_{1\ast}$ and $\sigma_1^{\ast}$, independent of $h$ and $E$, such that
 \begin{align}
 \sigma_{1\ast} a_1^E  (\c_h,\c_h)   \le a_{1h}^E(\c_h,\c_h) \le \sigma_1^{\ast} R_m^{-1}S_c\|\c_h\|_{1,E},\quad \forall \c_h \in \M(E).\label{CVEM-stability}
 \end{align}
\end{itemize}
\end{lemma}
\begin{proof}
It is clear that \refe{CVEM-cosistency} holds. For the stability, we need to prove that for all $\c_h \in  \M(E)$,
\begin{align*}
 \|(\I - \Pi_{k-1}^{0,E})\curl\,\c_h\|_{0,E} &\le \| \curl \, (\I - \Pi_k^{\nabla,E})\c_h \|_{0,E}\leq \sqrt{2}\| \nabla(\I - \Pi_k^{\nabla,E})\c_h \|_{0,E} ,\\
  \|(\I - \Pi_{k-1}^{0,E})\div\,\c_h\|_{0,E} &\le \| \div \, (\I - \Pi_k^{\nabla,E})\c_h \|_{0,E} \leq \sqrt{2}\| \nabla(\I - \Pi_k^{\nabla,E})\c_h \|_{0,E} ,
\end{align*}
where $\I$ is an identity operator. Making use of the properties of $L^2$-projection, we have
\begin{align*}
\|(\I - \Pi_{k-1}^{0,E})\curl\,\c_h\|_{0,E}^2 
& = ((\I - \Pi_{k-1}^{0,E})\curl\,\c_h,\curl\,(\I - \Pi_k^{\nabla,E})\c_h)_{E} \\
& \le \|(\I - \Pi_{k-1}^{0,E})\curl\,\c_h\|_{0,E}\|\curl\,(\I - \Pi_k^{\nabla,E})\c_h\|_{0,E},
\end{align*}
which is the desired result of the first inequality. Similarly, the second inequality can be derived. Thus, there hold
\begin{align*}
a_{1h}^E(\c_h,\c_h) & \le R_m^{-1} S_c\left(\|\Pi_{k-1}^{0,E}\curl\,\c_h\|_{0,E}^2 + \|\Pi_{k-1}^{0,E}\div\,\c_h\|_{0,E}^2 + \alpha_1^{\ast}\|\nabla (\I - \Pi_k^{\nabla,E})\c_h\|_{0,E}^2 \right)\\
& \le R_m^{-1}S_c(4 + \alpha_1^{\ast})\|\c_h\|_{1,E}^2 = R_m^{-1}S_c\sigma_1^{\ast}\|\c_h\|_{1,E}^2   \\
a_{1h}^E(\c_h,\c_h) & \geq  R_m^{-1}S_c \left(\|\Pi_{k-1}^{0,E}\curl\,\c_h\|_{0,E}^2 + \|\Pi_{k-1}^{0,E}\div\,\c_h\|_{0,E}^2 + \alpha_{1\ast}\|\nabla (\I - \Pi_k^{\nabla,E})\c_h\|_{0,E}^2 \right)\\
& \geq R_m^{-1}S_c \Big(\|\Pi_{k-1}^{0,E}\curl\,\c_h\|_{0,E}^2 + \|\Pi_{k-1}^{0,E}\div\,\c_h\|_{0,E}^2  + \frac{\alpha_{1\ast}}{4} \Big( \|(\I - \Pi_{k-1}^{0,E}) \curl\,\c_h\|_{0,E}^2  \\
& \quad +\|(\I - \Pi_{k-1}^{0,E})\div\,\c_h\|_{0,E}^2 \Big)\Big)  \\
& \geq R_m^{-1}S_c \min\{1,\frac{\alpha_{1\ast}}{4}\}\big( \|\curl\,\c_h\|_{0,E}^2 + \|\div\,\c_h\|_{0,E}^2  \big) = \sigma_{1\ast} a_1^E(\c_h,\c_h).
\end{align*}
The proof is completed.
\end{proof}
\subsection{The discrete formulation}
The virtual element formulation of the problem \refe{MHDrewriteeq}: $\forall (\v_h,\c_h,q_h) \in \V \times \M \times Q  $, find $(\u_h,\b_h,p_h) \in \V \times \M \times Q  $ such that
\begin{align}
\begin{cases}
A_h(\boldsymbol{\u}_h,\b_h;\v_h,\c_h) + B_h(\u_h,\b_h;\u_h,\b_h;\v_h,\c_h) - d_h(\v_h,\c_h;p_h)  = \left<F_h;\v_h,\c_h\right>,  \\
 d_h(\u_h,\b_h;q_h)  = 0,
\end{cases} \label{discrete-MHD}
\end{align}
where, $\forall \w_h \in \V, \ps_h \in \M$,
\begin{align*}
&A_h(\u_h,\b_h;\v_h,\c_h) = a_{0h}(\u_h,\v_h) + a_{1h}(\b_h,\c_h) = \sum_{E \in \T_h} a_{0h}^E(\u_h,\v_h) +  a_{1h}^E(\b_h,\c_h), \\
&B_h(\u_h,\b_h;\w_h,\ps_h;\v_h,\c_h) = a_{2h}(\u_h,\w_h,\v_h) - a_{3h}(\ps_h,\b_h,\v_h) + a_{3h}(\c_h,\b_h,\w_h), \\
&a_{2h} (\u_h,\w_h,\v_h) = \sum_{E \in \T_h} a_{2h}^E(\u_h,\w_h,\v_h)\\
& \qquad \qquad\qquad\;\;
=  \sum_{E\in \T_h} \frac{1}{2} \Big( (\Pi_{k-1}^{0,E}(\nabla \w_h)\Pi_k^{0,E}\u_h,\Pi_k^{0,E}\v_h)_{E}  - (\Pi_{k-1}^{0,E}(\nabla \v_h)\Pi_k^{0,E}\u_h,\Pi_k^{0,E}\w_h)_{E} \Big) ,\\
&a_{3h}(\ps_h,\b_h,\v_h) = \sum_{E\in \T_h} a_{3h}^E(\ps_h,\b_h,\v_h)=S_c \sum_{E\in \T_h}(\Pi_{k-1}^{0,E}\mathrm{curl}\, \ps_h \times \Pi_k^{0,E}\b_h ,\Pi_k^{0,E} \v_h)_{E} , \\
&d_{h}(\v_h,\c_h;q_h) = \sum_{E\in\T_h} d_{h}^E(\v_h,\c_h;q_h)= \sum_{E\in \T_h}(\div \, \v_h,q_h)_{E}  ,\\
&\left<F_h;\v_h ,\c_h\right> = \left(\f_h,\v_h\right) + (\g_h,\c_h) = \sum_{E\in\T_h} \left( (\Pi_k^{0,E}\f,\v_h )_{E} + (\Pi_k^{0,E}\g,\c_h)_{E}\right). 
\end{align*}
Obviously, the linear form $B_h(\cdot,\cdot;\cdot,\cdot;\cdot,\cdot)$ is antisymmetric, i.e. 
\begin{align*}
B_h(\cdot,\cdot;\w_h,\ps_h;\v_h,\c_h) = -B_h(\cdot,\cdot;\v_h,\c_h;\w_h,\ps_h), \quad \forall \w_h,\v_h\in \V,\ps_h,\c_h \in \M.
\end{align*} According to the definition of $\V$, it is easy to know that $\div_h\, \v_h \in Q$ for any $\v_h \in \V$, where $\div_h|_E = \div$. The discrete kernel space is defined as
\begin{align*}
 \Z_h = \{ \v_h \in \V: d_h(\v_h,\c_h;q_h) = 0 ,\quad \forall \,q_h \in Q\} =  \{ \v_h \in \V: \div_h \,\v_h = 0\}.
\end{align*}
Due to the discrete velocity $\u_h \in \Z_h$, then $\u_h$ is piecewise divergence-free. We define some norms
\begin{align*}
 |\v_h|_{1,h} = \sqrt{\sum_{E\in\T_h} |\v_h|_{1,E}}\, , \quad   \|(\v_h,\c_h)\|_{1,h} = \sqrt{ |\v_h|_{1,h}^2 + \|\c_h\|_1^2 }, \quad \forall \v_h \in \V, \c_h \in\M.
\end{align*}

There are some useful inequalities
\begin{align}
\|\Pi_k^{0,E} \v_h\|_{L^4(E)} \le \lmd_2 \| \v_h\|_{L^4(E) },&  \quad  \forall \v_h \in \V,\label{inquality1}\\
 \|\v_h\|_{L^4} \le  \lmd_3  |\v_h|_{1,h},&  \quad  \forall \v_h \in \H_0^1(\O)+ \V, \label{inquality2} \\
\|\v\|_{s-1,4,E} \le C \|\v\|_{s,E},&  \quad \forall \v \in \H^s(E), \label{inquality3}\\
\|(\I - \Pi_s^{0,E})\v\|_{L^4(E)} \le  C h_E^r |\v|_{r,4,E},& \quad \forall \v \in \W^{s+1,4}(E), \; 0\le r \le s+1,   \label{inquality4}
\end{align}
where $s$ is a positive integer. By using Inverse inequality of polynomials, the properties of $L^2$-projection and H$\ddot{\mathrm{o}}$lder inequality, we obtain \refe{inquality1}. The proof of \refe{inquality2}, see Lemma 0.1 of Appendix in \cite{Zhao-NS}. It is widely known that \refe{inquality3} holds. The inequality \refe{inquality4} can be found in \cite{DivFree-NS}. 

\begin{lemma}\label{dMHD-lemma} 
For any $\u_h,\w_h,\v_h \in \V,\,\b_h,\ps_h,\c_h \in \M$, there hold
\begin{align}
&|A_h(\u_h,\b_h;\v_h,\c_h)| \le \sigma^{\ast} \max\{R_{\nu}^{-1},R_m^{-1}S_c\} \|(\u_h,\boldsymbol{b}_h)\|_{1,h} \,\|(\v_h,\c_h)\|_{1,h} , \label{dMHDeq1}\\
&A_h(\u_h,\b_h;\u_h,\b_h) \geq \sigma_{\ast}\min\{ R_{\nu}^{-1},\lmd_1 R_m^{-1}S_c \} \|(\u_h,\b_h)\|_{1,h}^2,  \label{dMHDeq2}\\ 
&|B_h(\u_h,\b_h;\w_h,\ps_h;\v_h,\c_h)| \le  \hat{C} \|(\u_h,\b_h)\|_{1,h}\|(\w_h,\ps_h)\|_{1,h}\,\|(\v_h,\c_h)\|_{1,h}, \label{dMHDeq3} 
\end{align}
where $\sigma^{\ast} = \max\{\sigma_0^{\ast},\sigma_1^{\ast}\}, \sigma_{\ast}=\min\{\sigma_{0\ast},\sigma_{1\ast}\}$,  $\hat{C}$ is a constant independent of $h$.
\end{lemma}
\begin{proof}
By Cauchy-Schwarz inequality and stabilities \refe{NC-stability} and \refe{CVEM-stability}, we infer
\begin{align*}
|A_h(\u_h,\b_h;\v_h,\c_h)| &= |a_{0h}(\u_h,\v_h) + a_{1h}(\b_h,\c_h)| \\
&\le \sum_{E \in \T_h} \left( \sigma_0^{\ast} R_{\nu}^{-1} |\u_h|_{1,E}|\v_h|_{1,E} + \sigma_1^{\ast}R_m^{-1}S_c \|\b_h\|_{1,E}\|\c_h\|_{1,E} \right)\\
&\le \sigma^{\ast} \max\{R_{\nu}^{-1},R_m^{-1}S_c \} \|(\u_h,\b_h)\|_{1,h}\,\|(\v_h,\c_h)\|_{1,h}.
\end{align*}
Analogously, we get \refe{dMHDeq2}. Using \refe{inquality1} and \refe{inquality2} yields
\begin{align}
 a_{2h}(\u_h,\w_h,\v_h) &= \frac{1}{2} \sum_{E \in \T_h} \left( ((\Pi_{k-1}^{0,E} \nabla \w_h) \Pi_k^{0,E}\u_h,\Pi_k^{0,E}\v_h)_{E}  - ((\Pi_{k-1}^{0,E} \nabla \v_h) \Pi_k^{0,E}\u_h,\Pi_k^{0,E}\w_h)_{E}\right) \notag\\
& \le \frac{1}{2}\lmd_2^2 \sum_{E\in\T_h} \Big( \|\nabla \w_h\|_{0,E}\| \u_h\|_{L^4(E)} \|\v_h\|_{L^4(E)} + \|\nabla \v_h\|_{0,E}\| \u_h\|_{L^4(E)} \|\w_h\|_{L^4(E)} \Big) \notag\\
& \le \lmd_2^2\lmd_3^2  |\u_h|_{1,h}|\w_h|_{1,h}|\v_h|_{1,h}. \label{boundeq1}
\end{align} 
Similarly, we can obtain
\begin{align}
a_{3h}(\ps_h,\b_h,\v_h)  \le \sqrt{2} \lmd_0 \lmd_2^2\lmd_3 S_c  \|\ps_h\|_1 \|\b_h\|_1 |\v_h|_{1,h}. \label{boundeq2} 
\end{align}
Combining \refe{boundeq1}-\refe{boundeq2}, we derive \refe{dMHDeq3}. The proof is finished.
\end{proof}

\begin{lemma} \label{dinf-sup}
There exists a constant $\gm_1 > 0$ independent of $h$ such that
\begin{align}
\sup_{(\v_h,\c_h) \in \V\times\M} \frac{|d_h( \v_h,\c_h;q_h)|}{\|(\v_h,\c_h)\|_{1,h}} \geq \gm_1 \|q_h\|_0, \quad \forall q_h \in Q  . \label{dinf-supeq}
\end{align}
\end{lemma}
\begin{proof}
According to triangle inequality and Lemma \ref{NCVEM-interpolation}, it is not difficult to find that $|\v_I|_{1,h}  \le \eta |\v|_1$ for all $\v\in \H_0^1(\O)$, where $\eta$ is a positive constant independent of $h$. From the definition of interpolation (see \cite{Zhao-NS}) and the standard arguments in \cite{NS-1986}, for any $q_h \in Q$, we have
\begin{align*}
\sup_{(\v_I,\c_h) \in \V\times \M}  \frac{|d_h(\v_I,\c_h;q_h)|}{\|(\v_I,\c_h)\|_{1,h}} & \ge \sup_{(\v_I,\c_h) \in \V\times \boldsymbol{0}} \frac{|d_h(\v_I,\c_h;q_h)|}{\|(\v_I,\c_h)\|_{1,h}} \\
&= \sup_{\v_I \in \V} \frac{|(\div\,\v,q_h)|}{|\v_I|_{1,h}} \\
& \ge \frac{1}{\eta} \sup_{\v \in \H_0^1(\O)} \frac{|(\div\,\v,q_h)|}{|\v|_1}  \ge \frac{\gm_0}{\eta} \|q_h\|_0  = \gm_1 \|q_h\|_0.  
\end{align*}
\end{proof}
\begin{lemma}(\cite{FEM-2008}) \label{pi-interpolation}
Let $\v \in \H^l(E)$, $1\le l \le k+1$, there exists a polynomial $\v_{\pi} \in [\mathbb{P}_k(E)]^2$ such that 
\begin{align*}
\|\v - \v_{\pi} \|_{L^r(E)} + h|\v - \v_{\pi}|_{W^{1,r}(E)}  \le C h^{l}|\v|_{W^{l,r}(E)},\quad \forall 1\le r \leq \infty.
\end{align*} 
\end{lemma}
\begin{theorem}  \label{existeq}
Assuming 
\begin{align}
\frac{\hat{C}\|F_h\|_{\ast}}{(\sigma_{\ast}\min\{R_{\nu}^{-1}, \lmd_1 R_m^{-1}S_c\})^2} < 1,\label{existeq1}
\end{align} 
then the problem \refe{discrete-MHD} has a unique solution $(\u_h,\b_h,p_h) \in \V \times \M \times Q$ satisfying
\begin{align}
\sigma_{\ast}\min\{R_{\nu}^{-1}, \lmd_1 R_m^{-1}S_c\}\|(\u_h,\b_h)\|_{1,h} \le \|F_h\|_{\ast}. \label{existeq2}
\end{align}
\end{theorem}
\begin{proof}
Similar to \cite{G-MHD-1991,Dong-Three}, we can deduce \refe{existeq2} with just a few modifications.
\end{proof}
\setcounter{equation}{0}
\section{Error analysis}
\begin{lemma}\label{e-lemma1}
Under assumptions (\textbf{A1}) and (\textbf{A2}), suppose $\u\in \H_0^1(\O) \cap \H^{k+1}(\O)$, $\b\in \H_n^1(\O) \cap \H^{k+1}(\O)$, $\u_h \in \V$ and $\b_h \in \M$, for every $\v_h \in \V, \c_h \in \M$, there hold
\begin{align*}
&\Big| \sum_{E \in \T_h}a_2^E(\u,\u,\v_h) - a_{2h}(\u,\u,\v_h)\Big| \le Ch^k \|\u\|_{k+1}^2|\v_h|_{1,h},\\
& \Big|\sum_{E \in \T_h} a_{3}^E(\b,\b,\v_h) - a_{3h}(\b,\b,\v_h)\Big|   \le C h^k\|\b\|_k \|\b\|_{k+1} |\v_h|_{1,h}, \\
& \Big|\sum_{E \in \T_h} a_{3}^E(\c_h,\b,\u) - a_{3h}(\c_h,\b,\u)\Big|  \le C h^k\|\b\|_{k+1}\|\u\|_{k+1}\|\c_h\|_{1},\\
& \big|B_h(\u,\b;\u,\b;\v_h,\c_h) - B_h(\u_h,\b_h;\u_h,\b_h;\v_h,\c_h)\big|    \le  \hat{C} \Big(  \|(\u_h,\b_h)\|_{1,h} \|(\v_h,\c_h)\|_{1,h} \\
& \qquad \qquad \qquad \quad +  \|(\u - \u_h + \v_h,\b-\b_h + \c_h )\|_{1,h} \big(\|(\u,\b)\|_1  + \|(\u_h,\b_h)\|_{1,h} \big)    \Big) \|(\v_h,\c_h)\|_{1,h}.
\end{align*}
\end{lemma}
\begin{proof}
Similar to Lemma 4.3 in \cite{DivFree-NS}, with just a few modifications, the first inequality of Lemma \ref{e-lemma1} can be derived. According to the properties of $L^2$-projection, H$\ddot{\mathrm{o}}$lder inequality and \refe{inquality1}-\refe{inquality4}, we have
\begin{align*}
&\Big|\sum_{E \in \T_h}a_3^E(\b,\b,\v_h) - a_{3h}(\b,\b,\v_h)\Big|  \\
& = \Big|S_c\sum_{E\in\T_h} \int_E \Big( (\curl\,\b\times\b)\cdot \v_h  -  ( \Pi_{k-1}^{0,E}\curl\,\b\times \Pi_k^{0,E} \b)\cdot \Pi_k^{0,E}  \v_h\Big) \,\mathrm{d} E \Big|\\
& = \Big|S_c \sum_{E \in \T_h} \int_E \Big(  \big( (\I - \Pi_{k-1}^{0,E})(\curl\,\b\times\b) \big)\cdot (\I  - \Pi_k^{0,E})\v_h + \big( \curl\,\b\times(\I - \Pi_k^{0,E})\b \big)\cdot \Pi_k^{0,E}\v_h \\
&  \qquad + \big( (\I - \Pi_{k-1}^{0,E})\curl\,\b\times \Pi_k^{0,E}\b \big) \cdot \Pi_k^{0,E}\v_h \Big) \,\mathrm{d}E \Big| \\
& \le C h^k\sum_{E \in \T_h} \big(  \|\curl\,\b\times\b\|_{k-1,E}|\v_h|_{1,E} +  \|\b\|_{1,E}\|\b\|_{k,4,E}\|\v_h\|_{L^4(E)}  +  \|\b\|_{k+1,E} \|\b\|_{L^4(E)} \|\v_h\|_{L^4(E)}\big)\\
& \le Ch^k ( \|\b\|_1 + \|\b\|_k )\|\b\|_{k+1}|\v_h|_{1,h} \leq Ch^k \|\b\|_k \|\b\|_{k+1}|\v_h|_{1,h}.
\end{align*}
Similarly, we observe that
\begin{align*}
&\Big|\sum_{E \in \T_h}a_3^E(\c_h,\b,\u)  -  a_{3h}(\c_h,\b,\u)\Big|\\
 &= \Big|S_c \sum_{E \in \T_h} \int_E \Big((\I - \Pi_{k-1}^{0,E})\curl\,\c_h \cdot (\I - \Pi_{k-1}^{0,E})( \b\times \u) +  \Big( \Pi_{k-1}^{0,E} \curl\,\c_h \times (\I - \Pi_k^{0,E})\b \Big) \cdot \u \\
 & \quad + (\Pi_{k-1}^{0,E}\curl\,\c_h\times \Pi_k^{0,E} \b)\cdot (\I - \Pi_k^{0,E})\u \Big) \,\mathrm{d}E \Big|\\
& \le C h^k( \|\b\|_{k+1}\|\u\|_{k+1} + \|\b\|_{k+1}\|\u\|_1 + \|\b\|_1 \|\u\|_{k+1})\|\c_h\|_{1} \leq C h^k\|\b\|_{k+1}\|\u\|_{k+1}\|\c_h\|_{1}.
\end{align*}
By the antisymmetric of the linear form $B_h(\cdot,\cdot;\cdot,\cdot;\cdot,\cdot)$, we conclude that
\begin{align*}
& B_h(\u,\b;\u,\b;\v_h,\c_h) - B_h(\u_h,\b_h;\u_h,\b_h;\v_h,\c_h) \\
& = B_h(\u,\b;\u-\u_h,\b-\b_h;\v_h,\c_h) + B_h(\u - \u_h,\b - \b_h;\u_h,\b_h;\v_h,\c_h)  \\
& = B_h(\u,\b;\u-\u_h +\v_h,\b-\b_h +\c_h;\v_h,\c_h) \\
& \quad + B_h(\u - \u_h+\v_h,\b - \b_h+\c_h;\u_h,\b_h;\v_h,\c_h) -   B_h(\v_h,\c_h;\u_h,\b_h;\v_h,\c_h).
\end{align*}
Applying \refe{dMHDeq3}, the proof is completed. 
\end{proof}
\begin{lemma} \label{e-lemma2}
Under assumptions (\textbf{A1}) and (\textbf{A2}), let $\u \in \H_0^1(\O) \cap \H^{k+1}(\O), \w \in \H^2(\O)$ and $p \in L_0^2(\O) \cap H^k(\O)$, for every $\v_h \in \V$, we can get
\begin{align*}
|\mathcal{N}(\u,p;\v_h)| & \le Ch^k\left(R_{\nu}^{-1}\|\u\|_{k+1} + \|p\|_k + \|\u\|_{k+1}^2 \right)|\v_h|_{1,h}, \\
|\mathcal{N}(\u,p;\w,\v_h)| &\le Ch (R_{\nu}^{-1} \|\u\|_2 + \|p\|_1 + \|\u\|_2\|\w\|_2) |\v_h|_{1,h},
\end{align*}
where 
\begin{align*}
\mathcal{N}(\u,p;\v_h) &= \sum_{E \in \T_h} \int_{\partial E}\left(R_{\nu}^{-1}\nabla\u\cdot \n_E - p\n_E + \frac{1}{2}(\u\cdot \n_E)\u\right)\cdot\v_h \,\mathrm{d}\partial E, \\
\mathcal{N}(\u,p;\w,\v_h) &= \sum_{E \in \T_h} \int_{\partial E}  \left(R_{\nu}^{-1}   \nabla \u\cdot \n_E - p\n_E + \frac{1}{2} (\w\cdot \n_E) \u \right)\cdot \v_h  \,\mathrm{d}\partial E.
\end{align*}
\end{lemma}
\hspace*{-0.6cm} For the proof of Lemma \ref{e-lemma2}, see \cite{Liu-NS,Zhao-NS}.

\begin{theorem} \label{theorem-1norm}
Under assumptions (\textbf{A1}), (\textbf{A2}) and \refe{existeq1}, let $(\u,\b,p) \in ( \H_0^1(\O)\cap \H^{k+1}(\O)) \times (\H_n^1(\O)\cap \H^{k+1}(\O)) \times (L_0^2(\O) \cap H^k(\O))$ and $(\u_h,\b_h,p_h) \in \V\times \M \times Q$ be solutions of \refe{MHDrewriteeq} and \refe{discrete-MHD}, respectively. Then, the following estimates hold
\begin{align}
\|(\u-\u_h,\b -\b_h)\|_{1,h} &\le C(\u,\b,p,\f,\g,R_{\nu},R_m,S_c)h^k,  \label{theorem-u1normeq}\\ 
\|p - p_h\|_0 & \le C(\u,\b,p,\f,\g,R_{\nu},R_m,S_c)h^k.  \label{theorem-p0normeq} 
\end{align}

\begin{proof}
According to the definition of interpolation, which implies that $\u$ and $\u_I$ have the same degrees of freedom, we observe that if $\u \in \Z$, then $\u_I \in \Z_h$. Multiplying the first and second equations in \refe{MHDeq} by test functions $\delta_{\u} = \u_h - \u_I\in \Z_h$ and $\delta_{\b} =  \b_h - \b_I \in \M$, respectively, then using integration by parts, we deduce
\begin{align}
\sum_{E \in \T_h}A^E(\u,\b;\delta_{\u},\delta_{\b}) + \sum_{E \in \T_h}B^E(\u,\b;\u,\b;\delta_{\u},\delta_{\b}) - d_h(\delta_{\u},\delta_{\b};p)  = \left<F;\delta_{\u},\delta_{\b} \right>  + \mathcal{N}(\u,p;\delta_{\u}). \label{u1norm-eq1}
\end{align} 
Subtracting \refe{u1norm-eq1} from \refe{discrete-MHD}, taking $\v_h = \delta_{\u},\c_h = \delta_{\b}$ and using the fact $d_h(\delta_{\u},\delta_{\b};p)=0$, we derive
\begin{align}
\begin{aligned}
& A_h(\delta_{\u},\delta_{\b};\delta_{\u},\delta_{\b} )  \\
 &= \big\{A_h(\u_{\pi} - \u_I; \b_{\pi} - \b_I;\delta_{\u},\delta_{\b}) +  \sum_{E\in\T_h}A^E(\u - \u_{\pi},\b - \b_{\pi};\delta_{\u},\delta_{\b})\big\} \\
 &  \quad + \big\{ \sum_{E \in \T_h}B^E(\u,\b;\u,\b;\delta_{\u},\delta_{\b}) - B_h(\u,\b;\u,\b;\delta_{\u},\delta_{\b}) \big\}  + \big\{ B_h(\u,\b;\u,\b;\delta_{\u},\delta_{\b}) \\
 & \quad- B_h(\u_h,\b_h;\u_h,\b_h;\delta_{\u},\delta_{\b})\big\}   + \big\{\left< F_h;\delta_{\u},\delta_{\b}\right>  - \left<F;\delta_{\u},\delta_{\b}\right> \big\} - \mathcal{N}(\u,p;\delta_{\u})\\
 & = \Theta_1 + \Theta_2 + \Theta_3 + \Theta_4 + \Theta_5.
\end{aligned} \label{u1norm-eq2}
\end{align}
Applying triangle inequality, \refe{dMHDeq1}, Lemmas \ref{NCVEM-interpolation}, \ref{CVEM-interporlation} and \ref{pi-interpolation}, we obtain
\begin{align}
\begin{aligned}
\Theta_1
&= A_h(\u_{\pi} - \u_I,\b_{\pi} - \b_I;\delta_{\u},\delta_{\b}) + \sum_{E \in \T_h}A^E(\u - \u_{\pi},\b - \b_{\pi};\delta_{\u},\delta_{\b}) \\
& \le  C h^k\left( \|\u\|_{k+1} + \|\b\|_{k+1} \right)\|(\delta_{\u},\delta_{\b})\|_{1,h}.   
\end{aligned} \label{u1norm-eq3}
\end{align}
Using Lemma \ref{e-lemma1}, we have 
\begin{align}
\Theta_2
&= \sum_{E \in \T_h} B^E(\u,\b;\u,\b;\delta_{\u},\delta_{\b})- B_h(\u,\b;\u,\b;\delta_{\u},\delta_{\b}) \notag\\
& \le  Ch^k \big(\|\u\|_{k+1}^2 + \|\b\|_k\|\b\|_{k+1} + \|\u\|_{k+1}\|\b\|_{k+1} \big) \|(\delta_{\u},\delta_{\b})\|_{1,h}, \label{u1norm-eq4} \\
\Theta_3
&= B_h(\u,\b;\u,\b;\delta_{\u},\delta_{\b} ) - B_h(\u_h,\b_h;\u_h,\b_h;\delta_{\u},\delta_{\b})  \notag\\
& \le \hat{C} \|(\u_h,\b_h)\|_{1,h} \|(\delta_{\u},\delta_{\b})\|_{1,h}^2 + Ch^k \big(\|\u\|_{k+1} + \|\b\|_{k+1}\big)\big(\|(\u,\b)\|_1  + \|(\u_h,\b_h)\|_{1,h}\big) \|(\delta_{\u},\delta_{\b})\|_{1,h}   \notag\\ 
 & \le  \frac{\hat{C}  \|F_h\|_{\ast}}{\sigma_{\ast} \min\{R_{\nu}^{-1},\lmd_1  R_m^{-1}S_c\}} \|(\delta_{\u},\delta_{\b})\|_{1,h}^2 + C h^k \big( \|F\|_{\ast} + \|F_h\|_{\ast}  \big)  \big(\|\u\|_{k+1} + \|\b\|_{k+1}\big)\|(\delta_{\u},\delta_{\b})\|_{1,h}. \label{u1norm-eq5}
\end{align}
From the properties of $L^2$-projection, we know that
\begin{align}
\Theta_4
&= \left<F_h;\delta_{\u},\delta_{\b} \right>  - \left< F;\delta_{\u},\delta_{\b}\right>  \notag\\
& = \sum_{E \in \T_h} \int_E (\f - \Pi_k^{0,E} \f)\cdot (\I - \Pi_0^{0,E})\delta_{\u} + (\g - \Pi_k^{0,E} \g)\cdot (\I - \Pi_0^{0,E})\delta_{\b}\,\,\mathrm{d} E \notag\\
& \le C\sum_{E \in \T_h} \big( h_E^k |\f|_{k-1,E}|\delta_{\u}|_{1,E} + h_E^k |\g|_{k-1,E}|\delta_{\b}|_{1,E} \big)  \notag\\
& \le Ch^k \big( \|\f\|_{k-1} + \|\g\|_{k-1} \big)\|(\delta_{\u},\delta_{\b})\|_{1,h}.  
\label{u1norm-eq6}
\end{align}
Applying Lemma \ref{e-lemma2}, it is easy to see that
\begin{align}
\Theta_5 = - \mathcal{N}(\u,p;\delta_{\u}) \leq Ch^k \big(\|\u\|_{k+1} + \|p\|_k + \|\u\|_{k+1}^2\big)\|(\delta_{\u},\delta_{\b})\|_{1,h}. \label{u1norm-eq7}
\end{align}
Substituting \refe{u1norm-eq3}-\refe{u1norm-eq7} into \refe{u1norm-eq2} and using \refe{dMHDeq2}, we infer
\begin{align}
& \sigma_{\ast} \min\{R_{\nu}^{-1}, \lmd_1 R_m^{-1}S_c\} \Big( 1 - \frac{\hat{C}\|F_h\|_{\ast}}{(\sigma_{\ast}\min\{R_{\nu}^{-1},\lmd_1 R_m^{-1}S_c\})^2} \Big) \|(\delta_{\u},\delta_{\b})\|_{1,h} \notag\\
& \le C(\u,\b,p,\f,\g,R_{\nu},R_m,S_c)h^k .  \label{u1norm-eq8}
\end{align}
By \refe{u1norm-eq8}, assumption \refe{existeq1} and  triangle inequality, we can derive \refe{theorem-u1normeq}. 

Substacting \refe{u1norm-eq1} from \refe{discrete-MHD}, taking $\delta_{\u} = \v_h,\delta_{\b} = \c_h$, we obtain
\begin{align}
\begin{aligned}
&d_h(\v_h,\c_h;p) - d_h(\v_h,\c_h;p_h)\\
& = \Big\{\sum_{E\in\T_h}A^E(\u - \u_{\pi},\b - \b_{\pi} ; \v_h,\c_h) -  A_h(\u_h - \u_{\pi},\b_h - \b_{\pi};\v_h ,\c_h)\Big\} \\
& \quad +\Big\{ \sum_{E \in\T_h}B^E(\u,\b;\u,\b;\v_h,\c_h) - B_h(\u,\b;\u,\b;\v_h,\c_h)\Big\} + \big\{B_h(\u,\b;\u,\b;\v_h,\c_h) \\
& \quad- B_h(\u_h,\b_h;\u_h,\b_h;\v_h,\c_h)\big\}   + \big\{ \left<F_h;\v_h,\c_h \right> - \left<F;\v_h,\c_h\right> \big\} - \mathcal{N}(\u,p;\v_h) \\
&= \overline{\Theta}_1 + \overline{\Theta}_2+ \overline{\Theta}_3  +\overline{\Theta}_4+ \overline{\Theta}_5.
\end{aligned}\label{p0norm-eq1}
\end{align}
Similar to the estimates of the terms $\Theta_1,\Theta_2,\Theta_4,\Theta_5$, there hold
\begin{align}
\overline{\Theta}_1
& \le C \Big( \|(\u - \u_{h},\b - \b_{h})\|_{1,h}  + h^k\big(\|\u\|_{k+1} + \|\b\|_{k+1}\big)\Big) \|(\v_h,\c_h)\|_{1,h},\label{p0norm-eq2} \\
\overline{\Theta}_2
&\leq  Ch^k\Big(\|\u\|_{k+1}^2 + \|\b\|_k\|\b\|_{k+1} + \|\u\|_{k+1}\|\b\|_{k+1} \Big)\|(\v_h,\c_h)\|_{1,h}, \label{p0norm-eq3}\\
\overline{\Theta}_4 
& \leq Ch^k \big( \|\f\|_{k-1} + \|\g\|_{k-1} \big)\|(\v_h,\c_h)\|_{1,h},\label{p0nrom-eq4} \\
\overline{\Theta}_5 
& \leq Ch^k \big(\|\u\|_{k+1} + \|p\|_k + \|\u\|_{k+1}^2\big)\|(\v_h,\c_h)\|_{1,h}. \label{p0norm-eq5}
\end{align}
For the estimate of the term $\overline{\Theta}_3$, we observe that
\begin{align}
\begin{aligned}
\overline{\Theta}_3
&= B_h(\u,\b;\u,\b;\v_h,\c_h) - B_h(\u_h,\b_h;\u_h,\b_h;\v_h,\c_h)  \\
& = B_h(\u -\u_h,\b - \b_h;\u,\b;\v_h,\c_h) + B_h(\u_h,\b_h;\u - \u_h,\b - \b_h;\v_h,\c_h) \\
& \le C \|(\u - \u_h,\b - \b_h)\|_{1,h} \Big(\|(\u_h,\b_h)\|_{1,h} + \|(\u,\b)\|_{1}\Big)\|(\v_h,\c_h)\|_{1,h}.
\end{aligned}\label{p0norm-eq6}
\end{align}
Combining \refe{p0norm-eq2}-\refe{p0norm-eq6} and using \refe{theorem-u1normeq}, we conclude that
\begin{align}
& |d_h(\v_h,\c_h;p_h) - d_h(\v_h,\c_h;p)| \leq C(\u,\b,p,\f,\g,R_{\nu},R_m,S_c)h^k\|(\v_h,\c_h)\|_{1,h}. \label{p0norm-eq7}
\end{align}
Applying the discrete inf-sup condition \refe{dinf-supeq} yields
\begin{align*}
\|p_h -  p_I\|_0 &\le \frac{1}{\gm_1} \sup_{(\v_h,\c_h) \in \V \times \M} \frac{| d_h(\v_h,\c_h;p_h - p_I )|}{\|(\v_h,\c_h)\|_{1,h}}  \\
& =\frac{1}{\gm_1} \sup_{(\v_h,\c_h) \in \V \times \M} \frac{| d_h(\v_h,\c_h;p_h -  p )|}{\|(\v_h,\c_h)\|_{1,h}},
\end{align*}
which implies 
\begin{align*}
\|p_h - p_I\|_0 & \le C(\u,\b,p,\f,\g,R_{\nu},R_m,S_c)h^k.
\end{align*}
By triangle inequality and \refe{L2-eq}, we get \refe{theorem-p0normeq}. The proof is completed.
\end{proof}
\end{theorem}

For $L^2$-norm of error, we consider the following dual problem, see \cite{G-MHD-1991,NS-NFEM}.  
\begin{align}
\begin{cases}
- R_{\nu}^{-1} \Delta \w + \ds\frac{1}{2}\left(  (\nabla \u)^\mathrm{T}\w  - (\nabla \w)^\mathrm{T}\u \right) - (\nabla \w)\u + \nabla s + S_c \curl\,\Ph \times \b = \u - \u_h,  \quad \mathrm{in}\; \O, \\
 R_m^{-1}S_c \left( \curl\,\curl \,\Ph  - \nabla(\div\, \Ph)  \right) + S_c  \left(\curl\,\b \times \w - \curl\,\Ph \times \u - \curl (\b\times \w) \right) = \b  - \b_h,  \quad  \mathrm{in} \;\O ,\\
\div\,\w = 0 , \quad \mathrm{in} \;\O ,\\
\w = 0 ,\quad  \b \cdot \n = 0  ,\quad \curl\,\Ph \times \n  = 0,\quad \mathrm{on}\;\partial \O, 
\end{cases}  \label{dual-MHD}
\end{align}
where  $(\cdot)^\mathrm{T}$ represents the normal transpose. We assume that the solution of \refe{dual-MHD} satisfies $(\w,\Ph,s) \in \H^2(\O) \times \H^2(\O) \times H^1(\O)$ and there exists a constant $\lmd_4$ such that
\begin{align}
\|(\w,\Ph)\|_2 + \|s\|_1 \le \lmd_4 \|(\u-\u_h,\b - \b_h)\|_0.   \label{dual-eq}
\end{align}

\begin{theorem}  \label{theorem-0norm}
Under assumptions (\textbf{A1}), (\textbf{A2}) and \refe{existeq1}, let $(\u,\b,p) \in ( \H_0^1(\O) \cap \H^{k+1}(\O) )  \times ( \H_{n}^1(\O) \cap \H^{k+1}(\O) ) \times (L_0^2(\O) \cap H^k(\O))$ and $(\u_h,\b_h,p_h) \in \V \times \M \times Q$ be the solutions of \refe{MHDrewriteeq} and \refe{discrete-MHD}, respectively. If $(\w,\Ph,s) \in (\H^2(\O)\cap \H_0^1(\O)) \times (\H^2(\O)\cap \H_n^1(\O)) \times (H^1(\O)\cap L_0^2(\O))$, there holds
\begin{align*}
\|(\u - \u_h,\b - \b_h)\|_0 & \le C(\u,\b,p,\f,\g,R_{\nu},R_m,S_c)h^{k+1}.  
\end{align*}
\end{theorem}
\begin{proof}
Multiplying \refe{dual-MHD} by test functions $\u - \u_h$ and $\b - \b_h$, then adding \refe{discrete-MHD} and substracting \refe{u1norm-eq1}, with $\v_h=\delta_{\u} = \w_I \in \Z_h,\c_h=\delta_{\b}=\Ph_I$, we deduce 
\begin{align}
&\|(\u - \u_h,\b - \b_h)\|_0^2 \notag\\
&= \sum_{E \in \T_h}A^E(\u - \u_h,\b - \b_h;\w , \Ph) +  \sum_{E \in \T_h}B^E(\u,\b;\u-\u_h,\b - \b_h;\w,\Ph) \notag\\
& \quad+ \sum_{E \in \T_h}B^E(\u-\u_h,\b - \b_h;\u,\b;\w,\Ph)   -  \mathcal{N}(\w,s;\u,\u - \u_h) + A_h(\u_h,\b_h;\w_I,\Ph_I) \notag\\
& \quad + B_h(\u_h,\b_h;\u_h,\b_h;\w_I,\Ph_I) - \left< F_h;\w_I,\Ph_I\right> - \sum_{E \in \T_h} A^E(\u,\b;\w_I,\Ph_I) \notag\\
& \quad - \sum_{E \in \T_h} B^E(\u,\b;\u,\b;\w_I,\Ph_I) + \left<F;\w_I,\Ph_I \right> + \mathcal{N}(\u,p;\w_I) \notag\\
& = \sum_{E\in\T_h}A^E(\u-\u_h,\b -\b_h;\w - \w_I,\Ph - \Ph_I) + \Big\{A_h(\u_h,\b_h; \w_I ,\Ph_I) 
- \sum_{E\in \T_h}A^E(\u_h,\b_h; \w_I,\Ph_I) \Big\}   \notag\\
& \quad  + \Big\{B_h(\u_h,\b_h;\u_h,\b_h;\w_I - \w,\Ph_I - \Ph) - \sum_{E\in\T_h}B^E(\u,\b;\u,\b;\w_I - \w,\Ph_I - \Ph) \Big\}  \notag\\
& \quad+ \sum_{E\in\T_h}B^E(\u - \u_h,\b - \b_h;\u - \u_h,\b - \b_h;\w,\Ph) - \Big\{ \sum_{E\T_h}B^E(\u_h,\b_h;\u_h,\b_h;\w,\Ph)\notag\\
&\quad- B_h(\u_h,\b_h;\u_h,\b_h;\w,\Ph)\Big\} + \left<F - F_h; \w_I,\Ph_I \right>   + \big\{\mathcal{N}(\u,p;\w_I - \w) \notag\\
& \quad- \mathcal{N}(\w,s;\u,\u - \u_h)\big\} \notag\\
&= \Gamma_1 + \Gamma_2 + \Gamma_3 + \Gamma_4 + \Gamma_5 + \Gamma_6 + \Gamma_7 ,  \label{u0norm-eq1}
\end{align} 
where $\mathcal{N}(\u,p;\w)=0$. Applying Lemmas \ref{NCVEM-interpolation}-\ref{CVEM-interporlation} and \refe{dual-eq}, we get
\begin{align}
\Gamma_1 \le Ch\|(\u - \u_h,\b -\b_h )\|_{1,h} \|(\u - \u_h,\b - \b_h ) \|_0 . \label{u0norm-eq2}
\end{align}
Making use of consistencies \refe{NC-consistency}, \refe{CVEM-cosistency}, triangle inequality and \refe{dual-eq}, we obtain
\begin{align}
\Gamma_2
& = A_h(\u_h ,\b_h;\w_I, \Ph_I ) - \sum_{E\in \T_h}A^E(\u_h,\b_h; \w_I, \Ph_I) \notag\\
& = A_h(\u_h - \u_{\pi},\b_h - \b_{\pi};\w_I - \w_{\pi},\Ph_I - \Ph_{\pi} ) - \sum_{E\in\T_h}A^E(\u_h - \u_{\pi},\b_h - \b_{\pi}; \w_I - \w_{\pi}, \Ph_I - \Ph_{\pi})   \notag\\
& \le C \|(\u_h - \u_{\pi},\b_h - \b_{\pi} )\|_{1,h} \| (\w_I - \w_{\pi},\Ph_I - \Ph_{\pi})\|_{1,h} \notag \\
& \le  C \Big( h \|(\u - \u_h,\b - \b_h)\|_{1,h} + h^{k+1} \big(\|\u\|_{k+1} + \|\b\|_{k+1} \big)\Big)\|(\u - \u_h, \b - \b_h)\|_0 . \label{u0norm-eq3}
\end{align} 
From Lemmas \ref{dMHD-lemma}, \ref{e-lemma1}, \ref{NCVEM-interpolation}, \ref{CVEM-interporlation} and \refe{dual-eq}, we observe that
\begin{align}
\Gamma_3
&=B_h(\u_h,\b_h;\u_h,\b_h; \w_I - \w, \Ph_I - \Ph) - \sum_{E \in \T_h}B^E(\u,\b;\u,\b;\w_I - \w, \Ph_I - \Ph )  \notag\\
 & =B_h(\u_h - \u,\b_h - \b;\u_h,\b_h;\w_I - \w, \Ph_I - \Ph) + B_h(\u,\b;\u_h -\u,\b_h - \b; \w_I - \w, \Ph_I - \Ph)  \notag\\
 & \quad + B_h(\u,\b;\u,\b;\w_I - \w, \Ph_I - \Ph) - \sum_{E \in \T_h}B^E(\u,\b;\u,\b; \w - \w_I, \Ph_I - \Ph)  \notag\\
 & \le Ch \big(\|(\u,\b)\|_1 +  \|(\u_h,\b_h)\|_{1,h}  \big)\|(\u - \u_h,\b - \b_h)\|_{1,h}\|(\u - \u_h,\b - \b_h)\|_0\notag\\
 & \quad +C h^{k+1}\big( \|\u\|_{k+1}^2 + \|\b\|_{k}\|\b\|_{k+1} + \|\u\|_{k+1}\|\b\|_{k+1} \big) \|(\u - \u_h,\b - \b_h)\|_0. \label{u0norm-eq4}
\end{align}
Utilizing \refe{dMHDeq3}, Lemmas \ref{NCVEM-interpolation}-\ref{CVEM-interporlation} and \refe{dual-eq}, we derive
\begin{align}
\Gamma_4  \le C\|(\u - \u_h,\b - \b_h)\|_{1,h}^2 \|(\u - \u_h,\b - \b_h)\|_0.  \label{u0norm-eq5}
\end{align}
The term $\Gamma_5$ can be rewritten as
\begin{align}
\begin{aligned}
\Gamma_5 
&= -\sum_{E\in\T_h}B^E(\u_h,\b_h;\u_h,\b_h;\w,\Ph)+ B_h(\u_h,\b_h;\u_h,\b_h;\w,\Ph)\\
&= \Big\{a_{2h}(\u_h,\u_h,\w) - \sum_{E\T_h}a_2^E(\u_h,\u_h,\w)\Big\} + \Big\{
a_{3h}(\b_h,\b_h,\w)-\sum_{E\in\T_h}a_3^E(\b_h,\b_h,\w)   \Big\}\\
& \quad +\Big\{ a_{3h}(\Ph,\b_h,\u_h) -\sum_{E\in\T_h}a_3^E(\Ph,\b_h,\u_h) \Big\}\\
&= \Gamma_{51} + \Gamma_{52} + \Gamma_{53}.
\end{aligned} \label{u0norm-eq6}
\end{align}
Here, we only provide the estimates of $\Gamma_{52}$ and $\Gamma_{53}$, because the estimate of $\Gamma_{51}$ is analogous to $\Gamma_{52}$ and $\Gamma_{53}$. Similar to the proof of Theorem 4.7 in \cite{Liu-NS}, the term $\Gamma_{52}$ can be split as
\begin{align}
\Gamma_{52}
& = a_{3h}(\b_h,\b_h,\w ) - \sum_{E \in \T_h}a_3^E(\b_h,\b_h,\w)\notag\\
& = S_c \sum_{E \in \T_h} \Big( (\Pi_{k-1}^{0,E}\curl\,\b_h \times \Pi_k^{0,E}\b_h,\Pi_k^{0,E}\w)_{E} - (\curl\,\b_h\times \b_h,\w)_{E} \Big)\notag\\
& = S_c\sum_{E \in \T_h} \Big( ((\I - \Pi_{k-1}^{0,E})\curl\,\b_h \times (\b - \b_h),\w )_{E} +  (\Pi_{k-1}^{0,E}(\curl\,\b- \curl\,\b_h)\times(\I - \Pi_k^{0,E})\b_h,\w)_{E} \notag\\
&  \quad + (\Pi_{k-1}^{0,E}(\curl\,\b - \curl\,\b_h)\times \Pi_k^{0,E} \b_h,(\I - \Pi_k^{0,E})\w)_{E} -  (\Pi_{k-1}^{0,E}\curl\,\b\times(\I -\Pi_k^{0,E})\b_h,\w)_{E}\notag\\
&  \quad + (\Pi_{k-1}^{0,E}\curl\,\b\times \Pi_k^{0,E} (\b - \b_h),(\I - \Pi_k^{0,E})\w)_{E}  - ((\I - \Pi_{k-1}^{0,E} )\curl\,\b_h\times \b,\w)_{E} \notag\\
&  \quad - (\Pi_{k-1}^{0,E}\curl\, \b \times \Pi_k^{0,E}\b,(\I -\Pi_k^{0,E})\w)_{E}\Big)  \notag \\
& = S_c \left( \Gamma_{521} + \Gamma_{522} + \Gamma_{523} + \Gamma_{524} + \Gamma_{525} + \Gamma_{526} + \Gamma_{527}\right) . \label{u0norm-eq7-1}
\end{align}
From Lemma \ref{pi-interpolation}, \refe{Soboleveq} and \refe{dual-eq}, we get
\begin{align}
\begin{aligned}
\Gamma_{521} &= \sum_{E \in \T_h} ((\I - \Pi_{k-1}^{0,E})\curl\,\b_h\times (\b - \b_h),\w)_{E} \\
& \le \sum_{E \in \T_h} \|(\I - \Pi_{k-1}^{0,E})\curl\,(\b_h - \b + \b - \b_{\pi})\|_{0,E}\|\b - \b_h\|_{L^4(E)} \|\w\|_{L^4(E)} \\
& \le C\sum_{E \in \T_h} (|\b - \b_h|_{1,E} + |\b - \b_{\pi}|_{1,E}) \|\b - \b_h\|_{L^4(E)} \|\w\|_{L^4(E)} \\
& \le C(\|\b - \b_h\|_1 + h\|\b \|_2) \|\b - \b_h\|_1 \|(\u-\u_h,\b - \b_h)\|_0 . 
\end{aligned}\label{u0norm-eq7-2}
\end{align}
Using \refe{inquality4}, Lemma \ref{pi-interpolation}, \refe{Soboleveq} and \refe{dual-eq}, we arrive at
\begin{align}
\begin{aligned}
\Gamma_{522} &= \sum_{E \in \T_h} (\Pi_{k-1}^{0,E}(\curl\, \b - \curl\,\b_h)\times (\I - \Pi_k^{0,E})\b_h,\w)_{E} \\
& \le C\sum_{E \in \T_h}  \|\curl\,\b - \curl\,\b_h\|_{0,E} \|(\I - \Pi_k^{0,E})(\b_h - \b + \b - \b_{\pi})\|_{L^4(E)}  \|\w \|_{L^4(E)} \\
& \le C \sum_{E \in \T_h} |\b - \b_h|_{1,E} \|\b_h - \b + \b - \b_{\pi}\|_{L^4(E)}  \|\w \|_{L^4(E)} \\
& \le C (\|\b - \b_h\|_{1} + h\|\b\|_{2} ) \|\b - \b_h\|_{1}\|(\u-\u_h,\b - \b_h)\|_0 .
\end{aligned} \label{u0norm-eq7-3}
\end{align}
Obviously, it is easy to know that
\begin{align}
\Gamma_{523}  \le Ch  \|\b_h\|_1 \|\b - \b_h\|_1  \|(\u - \u_h,\b - \b_h)\|_0.  \label{u0norm-eq7-4}
\end{align}
Applying \refe{inquality1}, Lemma \ref{pi-interpolation}, \refe{Soboleveq} and \refe{dual-eq}, we deduce
\begin{align}
\begin{aligned}
\Gamma_{524} &=  -\sum_{E \in \T_h} (\Pi_{k-1}^{0,E}\curl\,\b\times (\I - \Pi_k^{0,E})\b_h,\w)_{E} \\
& \le C \sum_{E \in \T_h} \| \curl\, \b \|_{L^4(E)} \|(\I - \Pi_k^{0,E})(\b_h - \b + \b - \b_{\pi})\|_{0,E} \|\w\|_{L^4(E)}  \\
& \le  C \|\curl\, \b\|_{L^4}(h|\b - \b_h|_{1} + Ch^{k+1}|\b |_{k+1}) \|\w\|_{L^4}    \\
& \le C \|\b\|_2 \big( h\|\b - \b_h\|_1 + h^{k+1} \|\b\|_{k+1}\big)\|(\u - \u_h,\b - \b_h)\|_0.
\end{aligned}\label{u0norm-eq7-5}
\end{align}
Similar to \refe{u0norm-eq7-2}-\refe{u0norm-eq7-5}, we have
\begin{align}
\Gamma_{525}
& \le Ch \|\b\|_1 \|\b - \b_h\|_1 \|(\u - \u_h,\b - \b_h)\|_0,    \label{u0norm-eq7-6}\\
\Gamma_{526} 
&= -\sum_{E \in \T_h} ((\I - \Pi_{k-1}^{0,E})\curl\,\b_h\times \b,\w)_{E} \notag\\
& = -\sum_{E \in \T_h} ((\I - \Pi_{k-1}^{0,E})\curl\,\b_h,(\I - \Pi_0^{0,E})(\b\times \w))_{E} \notag\\
& \le C \sum_{E \in \T_h} h_E\|(\I - \Pi_{k-1}^{0,E})\curl\,(\b_h -\b + \b - \b_{\pi}) \|_{0,E} |\b\times \w|_{1,E} \notag\\
& \le C \|\b\|_2 \big(h  \|\b - \b_h\|_1  + h^{k+1} \|\b\|_{k+1} \big)\|(\u - \u_h,\b - \b_h)\|_0,     \label{u0norm-eq7-7}\\
\Gamma_{527} 
&= -\sum_{E \in \T_h} (\Pi_{k-1}^{0,E}\curl\,\b\times \Pi_k^{0,E}\b,(\I - \Pi_k^{0,E})\w)_{E} \notag\\
& =  -\sum_{E \in \T_h} ((\I - \Pi_{k-1}^{0,E})(\Pi_{k-1}^{0,E}\curl\,\b \times \Pi_k^{0,E}\b),(\I - \Pi_k^{0,E}) \w)_{E} \notag \\
& \le C\sum _{E \in \T_h} h_E^{k+1}  |\Pi_{k-1}^{0,E} \curl\,\b \times \Pi_k^{0,E} \b |_{k-1,E} |\w|_{2,E} \notag\\
& \le Ch^{k+1} \|\b\|_{k} \|\b\|_{k+1} \|(\u - \u_h, \b - \b_h)\|_0. \label{u0norm-eq7-8}
\end{align}
Combining \refe{u0norm-eq7-2}-\refe{u0norm-eq7-8}, we derive
\begin{align}
\Gamma_{52} \leq C\Big( h^{k+1} + h\|(\u-\u_h,\b-\b_h)\|_{1,h} + \|(\u-\u_h,\b-\b_h)\|_{1,h}^2 \Big)\|(\u-\u_h,\b-\b_h)\|_0. \label{u0norm-eq7-9}
\end{align}
Analogously, by adding and subtracting terms, the term $\Gamma_{53}$ can be written as
\begin{align}
\Gamma_{53}
&= a_{3h}(\Ph,\b_h,\u_h) - \sum_{E \in \T_h}a_3^E(\Ph,\b_h,\u_h) \notag\\
& = S_c \sum_{E \in \T_h} \Big( (\Pi_{k-1}^{0,E}\curl\,\Ph\times \Pi_k^{0,E} \b_h,\Pi_k^{0,E}  \u_h)_{E} - (\curl\,\Ph\times\b_h,\u_h)_{E}\Big) \notag\\
& = S_c\sum_{E \in \T_h} \Big( -(\Pi_{k-1}^{0,E}\curl\,\Ph\times(\I - \Pi_{k}^{0,E})\b_h,\u_h)_{E} - (\Pi_{k-1}^{0,E} \curl\,\Ph\times \Pi_k^{0,E} \b_h,(\I - \Pi_k
^{0,E})\u_h)_{E} \notag\\
&  \quad + ((\I - \Pi_{k-1}^{0,E})\curl\,\Ph\times(\b - \b_h),\u_h)_{E} +  ((\I - \Pi_{k-1}^{0,E})\curl\,\Ph \times \b,\u - \u_h)_{E} \notag\\
& \quad - ((\I - \Pi_{k-1}^{0,E})\curl\,\Ph\times \b,\u)_{E}\Big)\notag\\
& = S_c\left( \Gamma_{531} + \Gamma_{532} + \Gamma_{533} +\Gamma_{534} +\Gamma_{535}\right) . \label{u0norm-eq8-1}
\end{align}
Similar to \refe{u0norm-eq7-2}-\refe{u0norm-eq7-8}, we deduce
\begin{align}
\begin{aligned}
\Gamma_{531} 
&  \le C|\u_h|_{1,h}\big(h\|\b - \b_h\|_1 + h^{k+1} \|\b\|_{k+1}\big) \|(\u - \u_h,\b - \b_h)\|_0,  \\
\Gamma_{532}
&\le C \|\b_h\|_1 \big( h |\u - \u_h|_{1,h}+ h^{k+1} \|\u\|_{k+1}\big)\|(\u - \u_h,\b - \b_h)\|_0 , \\
\Gamma_{533} 
& \le C h |\u_h|_{1,h} \|\b - \b_h\|_1 \|(\u - \u_h,\b - \b_h)\|_0, \\
\Gamma_{534} 
& \le C h\|\b\|_1 |\u - \u_h|_{1,h} \|(\u - \u_h,\b - \b_h)\|_0, \\
\Gamma_{535} 
& \le Ch^{k+1}\|\u\|_{k+1} \|\b\|_{k+1} \|(\u - \u_h,\b - \b_h)\|_0 . 
\end{aligned} \label{u0norm-eq8-2}
\end{align}
Substituting \refe{u0norm-eq8-2} into \refe{u0norm-eq8-1} yields
\begin{align}
\Gamma_{53} 
& \le C\Big(h^{k+1} + h\|(\u - \u_h,\b - \b_h)\|_{1,h} \Big)\|(\u - \u_h, \b - \b_h)\|_0.  \label{u0norm-eq8-3}
\end{align}
By the same way, following estimate can be derived 
\begin{align}
\Gamma_{51} \leq C \Big(h^{k+1} +  h\|(\u - \u_h,\b - \b_h)\|_{1,h} + \|(\u - \u_h,\b - \b_h)\|_{1,h}^2\Big)\|(\u - \u_h,\b - \b_h)\|_0.   \label{u0norm-eq9}
\end{align}
According to the properties of $L^2$-projection and triangle inequality, we find 
\begin{align}
\begin{aligned}
\Gamma_6
&= \left<F - F_h; \w_I, \Ph_I\right>  \\
& =  \sum_{E \in \T_h} \left( (\f - \Pi_k^{0,E}\f,\w_I - \w_{\pi})_{E} + (\g - \Pi_k^{0,E}\g,\Ph_I - \Ph_{\pi})_{E} \right)  \\
& \le C\sum_{E \in \T_h} h_E^{k+1} \left( |\f|_{k-1,E} |\w|_{2,E} + |\g|_{k-1,E} |\Ph|_{2,E}  \right)\\
& \le Ch^{k+1} \left(\|\f\|_{k-1} + \|\g\|_{k-1}  \right)\|(\u - \u_h, \b - \b_h)\|_0.
\end{aligned}\label{u0norm-eq10}
\end{align}
Using Lemmas \ref{e-lemma2}, \ref{NCVEM-interpolation} and \refe{dual-eq}, the estimate of the term $\Gamma_7$ is as follows:
\begin{align}
\begin{aligned}
\Gamma_7 &\le Ch^{k+1}\left(\|\u\|_{k+1} + \|p\|_k + \|\u\|_{k+1}^2 \right)\|(\u - \u_h,\b - \b_h)\|_0  \\
& \quad + Ch(1 + \|\u\|_2)|\u - \u_h|_{1,h} \|(\u - \u_h,\b - \b_h)\|_0.
\end{aligned}   \label{u0norm-eq11}
\end{align}
Combining \refe{u0norm-eq2}-\refe{u0norm-eq5}, \refe{u0norm-eq7-9} and \refe{u0norm-eq8-3}-\refe{u0norm-eq11} and using Theorem \ref{theorem-u1normeq}, we obtain 
\begin{align*}
\|(\u - \u_h,\b - \b_h)\|_0^2 \leq C(\u,\b,p,\f,\g,R_{\nu},R_m,S_c)h^{k+1}. 
\end{align*} 
The proof is finished.  
\end{proof}

\section{Numerical examples}
In this section, we first introduce a fixed-point algorithm (Oseen iteration) for the handing of nonlinear problems. Next, by the Example \ref{example1}, we verify the convergence rates of the proposed nonconforming VEM for the cases $k=1,2$, and compare the nonconforming VEM, the conforming VEM and the FEM for the case $k=1$. The above conforming VEM only changes the approximation space of the velocity compared with the nonconforming VEM, the velocity is approximate by the low-order Stokes-type virtual element \cite{VEM-MHD-t-2D}. For the FEM, the stable finite element pair $(P_{1b},P_1,P_1)$ is used for discretization. For both the VEM and FEM, we adopt Oseen iteration  to deal with nonlinear system.  The Oseen iteration scheme of the VEM  is as follows:
\begin{align*}
\begin{cases}
A_h(\boldsymbol{\u}_h^n,\b_h^n;\v_h,\c_h) + B_h(\u_h^{n-1},\b_h^{n-1};\u_h^n,\b_h^n;\v_h,\c_h) - d_h(\v_h,\c_h;p_h^n)  = \left<F_h;\v_h,\c_h\right>,  \\
 d_h(\u_h^n,\b_h^n;q_h)  = 0,
\end{cases} 
\end{align*}
with the iterative initial value $\u_h^0=\mathbf{0},\b_h^0 =\mathbf{0}$, the iterative tolerance is 1e-7.  The Oseen iteration of the FEM can be found in \cite{Dong-Three}. Then, in Example \ref{example2}, we investigate a benchmark problem. All experiments are implemented with FEALPy package \cite{Wei-Fealpy}.

In order to compute VEM errors, we consider the following computable error quantities:
\begin{align*}
\|e_{\u}\|_0  &= \sqrt{\sum_{E \in \T_h} \|\u - \Pi_k^{0,E} \u_h\|_{0,E}^2},  \quad |e_{\u}|_1 = \sqrt{ \sum_{E \in \T_h} |\u - \Pi_k^{\nabla,E} \u_h|_{1,E}^2},\quad \|e_p\|_0 = \|p - p_h \|_0,\\
 \|e_{\b}\|_0  &= \sqrt{\sum_{E \in \T_h} \|\b - \Pi_k^{0,E} \b_h\|_{0,E}^2},  \quad \|e_{\b}\|_1 = \sqrt{ \sum_{E \in \T_h} \big(\|\b - \Pi_k^{0,E}\b_h\|_{0,E}^2 + |\b - \Pi_k^{\nabla,E} \b_h |_{1,E}^2\big) }.
\end{align*}
For the choice of the stabilization terms $\S_0^E,\S_1^E$, we follow \cite{Zhao-NS,Equivalent-VEM}
\begin{align*}
\S_0^E(\u_h,\v_h) = \sum_{i=1}^{\mathrm{dim}(\V(E))} \mathrm{dof}_i^{\V(E)}(\u_h)\mathrm{dof}_i^{\V(E)}(\v_h),\\
\S_1^E(\b_h,\c_h) = \sum_{i=1}^{\mathrm{dim}(\M(E))} \mathrm{dof}_i^{\M(E)}(\b_h)\mathrm{dof}_i^{\M(E)}(\c_h),
\end{align*}
where $\mathrm{dof}_i^{\W}(\v)$ is the $i$th degree of freedom of smooth enough function $\v$ to the space $\W$, here $\W$ is $\V(E)$ or $\M(E)$.
Of course, the choice is not unique, there are also some other options \cite{Ill-VEM,Small-edge,Chen-VEM}.
\begin{figure}[H]
\vspace{-0.5cm}
\centering
\subfigure
{
    \begin{minipage}[b]{.3\linewidth}
        \centering
        \includegraphics[scale=0.35]{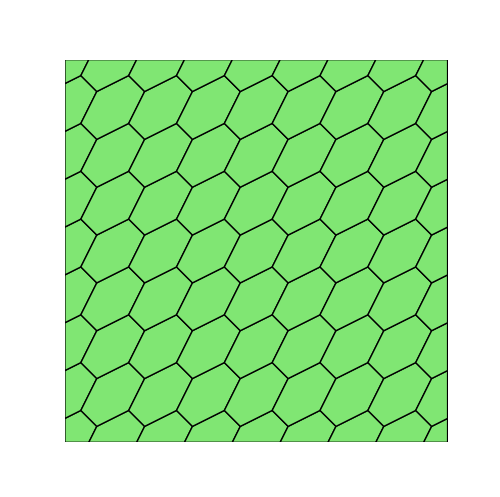}
    \end{minipage}
}
\subfigure
{
 	\begin{minipage}[b]{.3\linewidth}
        \centering
        \includegraphics[scale=0.35]{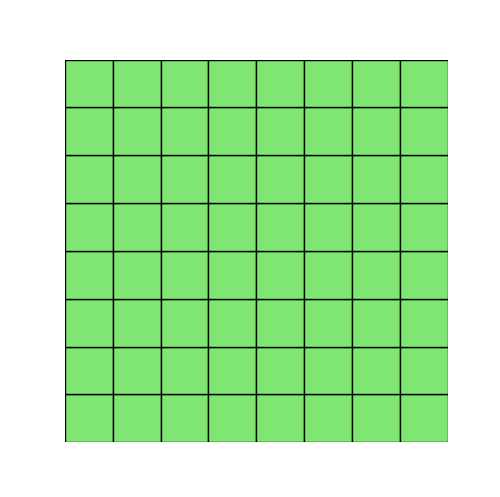}
    \end{minipage}
}
\subfigure
{
 	\begin{minipage}[b]{.3\linewidth}
        \centering
        \includegraphics[scale=0.35]{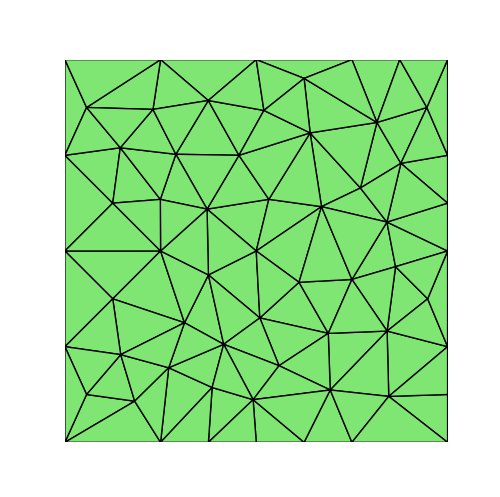}
    \end{minipage}
}
\vspace{-0.3cm}
\caption{Sample meshes: $\T_h^1$(left), $\T_h^2$(middle) and $\T_h^3$(right).} \label{Fig1}
\end{figure}

\begin{example} \label{example1}
We consider solving MHD with $R_{\nu}=R_m=S_c=1$ on domain $\O = [0,1]^2$. The source terms $\f$ and $ \g$ satisfy the exact solutions
\begin{align*}
& \u = \begin{pmatrix}
 \sin(\pi x_1)^2 \sin(\pi x_2) \cos(\pi x_2) \\
- \sin(\pi x_1) \sin(\pi x_2)^2 \cos(\pi x_1)
\end{pmatrix} ,\; \b = \begin{pmatrix}
 \sin(\pi x_1) \cos(\pi x_2) \\
-\sin(\pi x_2)  \cos(\pi x_1)
\end{pmatrix}, \; p = \cos(\pi x_1) \cos(\pi x_2) .
\end{align*}

In Table \ref{table1}, we list the errors and the optimal convergence rates of the velocity on meshes $\T_h^1$ and $\T_h^2$ (see Figure \ref{Fig1}) for the cases $k=1,2$, as well as $L^2$-norm of $\div_h \u_h$. We observe that the discrete velocity is divergence-free up to machine precision. Table \ref{table2} displays the errors and the optimal convergence rates of the magnetic and pressure on meshes $\T_h^1$ and $\T_h^2$ for the cases $k=1,2$. Tables \ref{table1} and \ref{table2} illustrate validity of the theoretical analysis. 
\begin{table}[H] 
\vspace*{-0.3cm}
\centering
\caption{Convergence rates of velocity on meshes $\T_h^1$ and $\T_h^2$.}
\begin{tabular}{cccclclcl}
\hline
meshes & $k$ & $h$  & $\|e_{\u}\|_0$ & Rate & $|e_{\u}|_1$ & Rate   & $\|\div_h \u_h\|_0$ \\ \hline 
\multirow{5}*{$\T_h^1$} & \multirow{5}*{1} 
  & 0.3727 & 7.9642e-02 & -    & 1.4316e+00 & -    & 1.8392e-15    \\
& & 0.1863 & 2.7525e-02 & 1.53 & 7.0057e-01 & 1.03 & 3.6925e-15     \\
& & 0.0932 & 7.8488e-03 & 1.81 & 3.3440e-01 & 1.07 & 2.7869e-15     \\
& & 0.0466 & 2.0667e-03 & 1.93 & 1.6374e-01 & 1.03 & 1.1629e-14     \\
& & 0.0233 & 5.2739e-04 & 1.97 & 8.1386e-02 & 1.01 & 4.7360e-14 \\ \hline
\multirow{5}*{$\T_h^1$} & \multirow{5}*{2} 
  & 0.3727 & 2.8702e-02 & -    & 1.4227e+00 & -    & 6.2575e-15 \\ 
& & 0.1863 & 3.7358e-03 & 2.94 & 5.5399e-01 & 1.36 & 1.1871e-14 \\   
& & 0.0932 & 4.2119e-04 & 3.15 & 1.5298e-01 & 1.86 & 2.5457e-14 \\    
& & 0.0466 & 4.4875e-05 & 3.23 & 3.9484e-02 & 1.95 & 4.4265e-14 \\    
& & 0.0233 & 4.9710e-06 & 3.17 & 9.9924e-03 & 1.98 & 1.2225e-13 \\ \hline 
\multirow{5}*{$\T_h^2$} & \multirow{5}*{1} 
  & 0.3536 & 1.5356e-01 & -    & 1.6980e+00 & -    & 1.0977e-15  \\ 
& & 0.1768 & 4.4590e-02 & 1.78 & 7.6794e-01 & 1.14 & 3.0290e-15  \\   
& & 0.0884 & 1.1317e-02 & 1.98 & 3.3497e-01 & 1.2  & 3.9649e-15  \\   
& & 0.0442 & 2.8199e-03 & 2.0  & 1.5790e-01 & 1.09 & 8.4733e-15  \\ 
& & 0.0221 & 7.0356e-04 & 2.0  & 7.7581e-02 & 1.03 & 3.3386e-14  \\ \hline
\multirow{5}*{$\T_h^2$} & \multirow{5}*{2} 
  & 0.3536 & 5.6481e-02 & -    & 1.5316e+00 & -    & 8.3672e-15  \\  
& & 0.1768 & 9.0257e-03 & 2.65 & 8.8572e-01 & 0.79 & 8.7864e-15  \\   
& & 0.0884 & 9.1090e-04 & 3.31 & 3.0408e-01 & 1.54 & 1.5897e-14  \\  
& & 0.0442 & 7.7982e-05 & 3.55 & 8.2445e-02 & 1.88 & 4.1011e-14  \\   
& & 0.0221 & 6.5169e-06 & 3.58 & 2.1539e-02 & 1.94 & 1.1272e-13 \\ \hline 
\end{tabular}  \label{table1}
\end{table} 
\begin{table}[H] 
\centering
\vspace*{-0.3cm}
\caption{Convergence rates of magnetic field and pressure on meshes $\T_h^1$ and $\T_h^2$.}
\begin{tabular}{cccclclclc}
\hline
meshes & $k$ & $h$  & $\|e_{\b}\|_0$ & Rate & $\|e_{\b}\|_1$ & Rate   & $\|e_p\|_0$ & Rate\\ \hline 
\multirow{5}*{$\T_h^1$} & \multirow{5}*{1} 
  & 0.3727 & 4.7679e-02 & -    & 9.1759e-01 & -    & 3.5180e-01 & -      \\
& & 0.1863 & 1.3848e-02 & 1.78 & 4.9783e-01 & 0.88 & 1.7097e-01 & 1.04   \\
& & 0.0932 & 3.7189e-03 & 1.9  & 2.5769e-01 & 0.95 & 6.4053e-02 & 1.42   \\
& & 0.0466 & 9.5814e-04 & 1.96 & 1.3086e-01 & 0.98 & 2.5234e-02 & 1.34   \\
& & 0.0233 & 2.4266e-04 & 1.98 & 6.5904e-02 & 0.99 & 1.1266e-02 & 1.16   \\ \hline
\multirow{5}*{$\T_h^1$} & \multirow{5}*{2} 
  & 0.3727 & 7.6370e-03 & -    & 1.7843e-01 & -    & 1.8498e-01 & -      \\ 
& & 0.1863 & 1.0430e-03 & 2.87 & 4.8383e-02 & 1.88 & 4.9715e-02 & 1.9    \\   
& & 0.0932 & 1.3037e-04 & 3.0  & 1.2346e-02 & 1.97 & 9.4398e-03 & 2.4    \\    
& & 0.0466 & 1.6216e-05 & 3.01 & 3.1185e-03 & 1.99 & 1.5328e-03 & 2.62   \\    
& & 0.0233 & 2.0260e-06 & 3.0  & 7.8440e-04 & 1.99 & 2.5201e-04 & 2.6     \\ \hline 
\multirow{5}*{$\T_h^2$} & \multirow{5}*{1} 
  & 0.3536 & 7.8227e-02 & -    & 9.9789e-01 & -    & 3.1471e-01 & -     \\ 
& & 0.1768 & 1.9856e-02 & 1.98 & 5.0244e-01 & 0.99 & 1.5002e-01 & 1.07  \\  
& & 0.0884 & 4.9809e-03 & 2.0  & 2.5167e-01 & 1.0  & 5.6183e-02 & 1.42  \\   
& & 0.0442 & 1.2462e-03 & 2.0  & 1.2589e-01 & 1.0  & 2.2761e-02 & 1.3   \\ 
& & 0.0221 & 3.1162e-04 & 2.0  & 6.2955e-02 & 1.0  & 1.0399e-02 & 1.13  \\ \hline
\multirow{5}*{$\T_h^2$} & \multirow{5}*{2} 
  & 0.3536 & 6.3629e-03 & -    & 1.9293e-01 & -    & 7.4471e-01 & -     \\  
& & 0.1768 & 8.0376e-04 & 2.98 & 4.8794e-02 & 1.98 & 2.0652e-01 & 1.85  \\   
& & 0.0884 & 9.6233e-05 & 3.06 & 1.2098e-02 & 2.01 & 3.4511e-02 & 2.58  \\  
& & 0.0442 & 1.1679e-05 & 3.04 & 3.0013e-03 & 2.01 & 5.1169e-03 & 2.75  \\   
& & 0.0221 & 1.4478e-06 & 3.01 & 7.4778e-04 & 2.0  & 9.0330e-04 & 2.5   \\ \hline 
\end{tabular}  \label{table2}
\end{table} 
\begin{figure}[H]
\vspace{-0.5cm}
\centering
\subfigure
{
    \begin{minipage}[b]{.4\linewidth}
        \centering
        \includegraphics[scale=0.4]{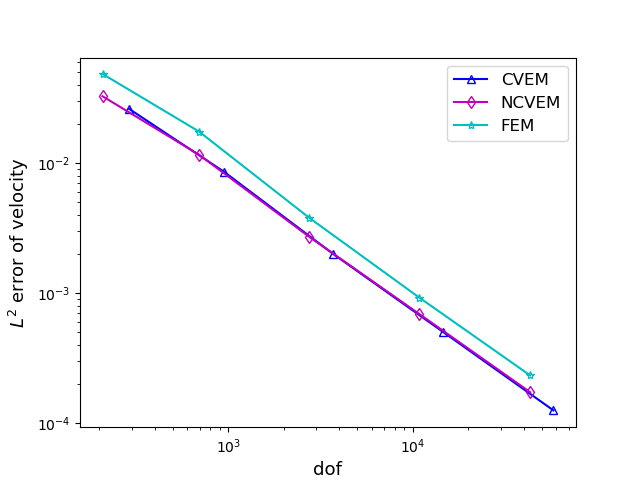}
    \end{minipage}
}
\subfigure
{
 	\begin{minipage}[b]{.4\linewidth}
        \centering
        \includegraphics[scale=0.4]{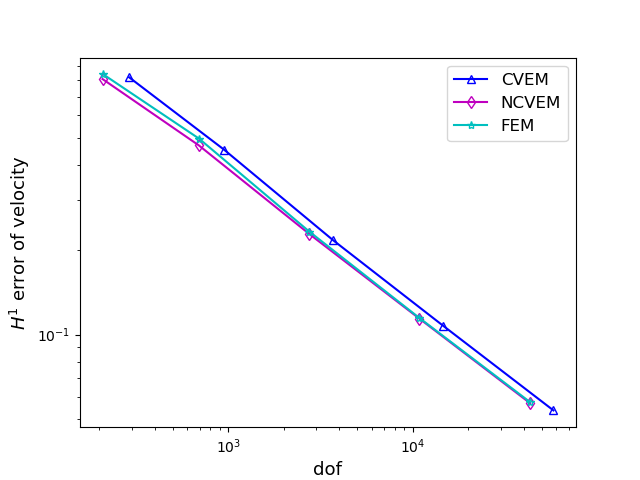}
    \end{minipage}
}
\subfigure
{
    \begin{minipage}[b]{.4\linewidth}
        \centering
        \includegraphics[scale=0.4]{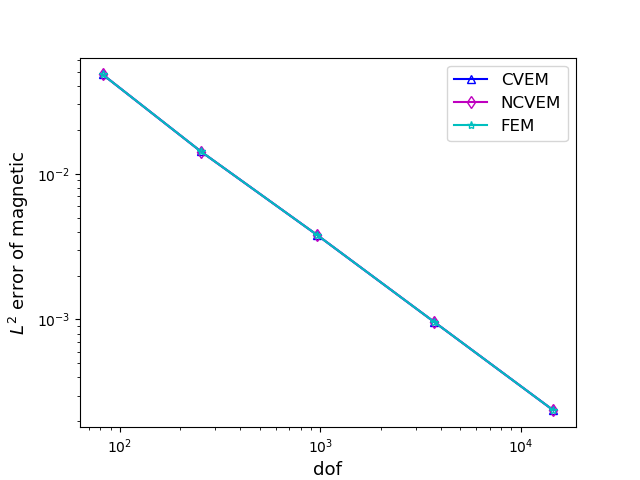}
    \end{minipage}
}
\subfigure
{
 	\begin{minipage}[b]{.4\linewidth}
        \centering
        \includegraphics[scale=0.4]{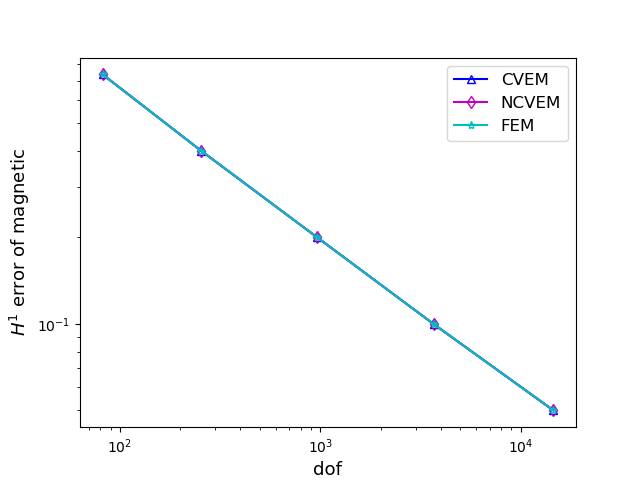}
    \end{minipage}
}
\caption{Errors of velocity and magnetic field on meshes $\T_h^3$.} \label{Fig2}
\end{figure}
\begin{figure}[H]
\vspace{-0.5cm}
\centering
\subfigure
{
    \begin{minipage}[b]{.4\linewidth}
        \centering
        \includegraphics[scale=0.4]{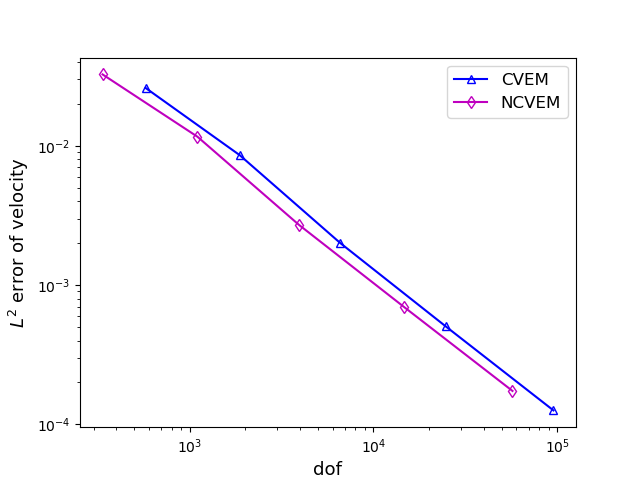}
    \end{minipage}
}
\subfigure
{
 	\begin{minipage}[b]{.4\linewidth}
        \centering
        \includegraphics[scale=0.4]{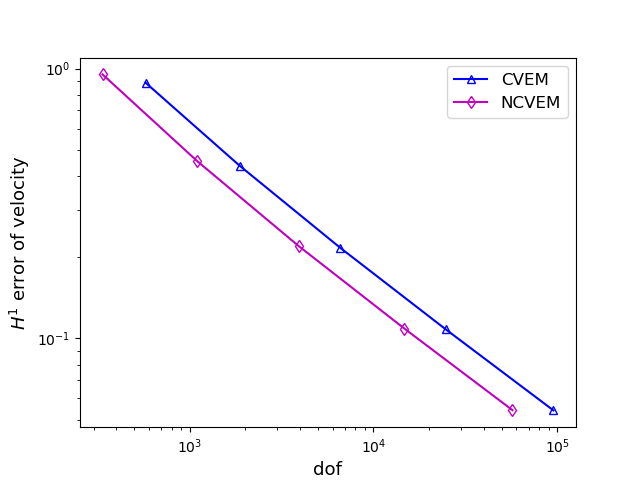}
    \end{minipage}
}
\subfigure
{
    \begin{minipage}[b]{.4\linewidth}
        \centering
        \includegraphics[scale=0.4]{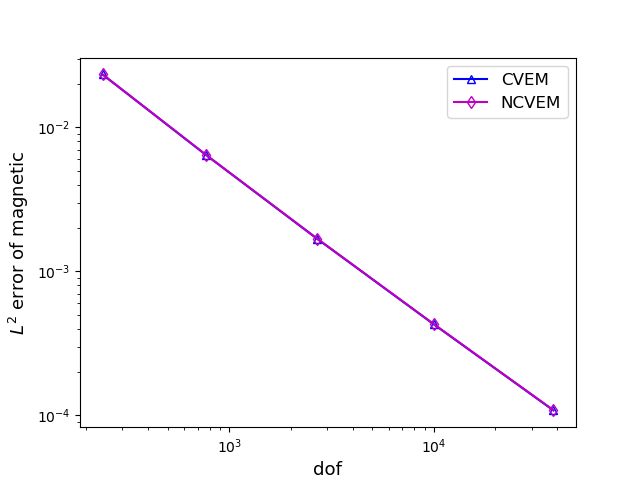}
    \end{minipage}
}
\subfigure
{
 	\begin{minipage}[b]{.4\linewidth}
        \centering
        \includegraphics[scale=0.4]{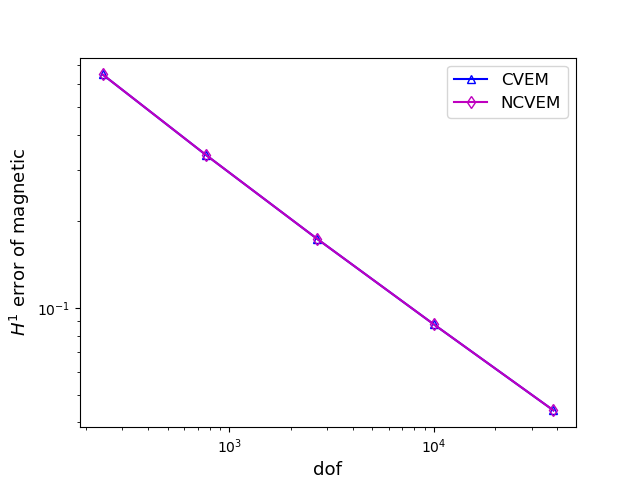}
    \end{minipage}
}
\caption{Errors of velocity and magnetic field on meshes $\T_h^1$.} \label{Fig3}
\end{figure}
\end{example}

For comparison, we also solve the MHD problem by the conforming VEM and the FEM.
In Figure \ref{Fig2}, we compare with three methods on meshes $\T_h^3$ for the case $k=1$, by displaying the error curves of the velocity and magnetic field with respect to the number of degrees of freedom. Figure \ref{Fig3} shows the compared results of the conforming VEM and the nonconforming VEM on meshes $\T_h^1$ for the case $k=1$. From Figures \ref{Fig2} and \ref{Fig3}, we observe that our proposed nonconforming VEM can provide better approximations for the velocity than the conforming VEM and the FEM in the numerical test for the case $k=1$.

\begin{example} \label{example2}
We consider Hartmann flow in the channel $\O= [0,6]\times[-1,1]$. The transverse magnetic field $\b_d=(0,1)$ is applied. The solutions of MHD are 
\begin{align*}
\u(x_1,x_2) = (u(x_2),0) ,\quad \b(x_1,x_2) = (b(x_2),1), \quad p(x_1,x_2) = -Gx_1 - S_cb^2(x_2)/2,
\end{align*}
where
\begin{align*}
u(x_2) = \frac{GR_{\nu}}{\Ha\cdot \tanh(\Ha)}\left( 1 - \frac{\cosh(x_2 \Ha)}{\cosh(\Ha)} \right), \quad b(x_2) = \frac{G}{S_c}\left(\frac{\sinh(x_2 \Ha)}{\sinh(\Ha)} - x_2  \right).
\end{align*}
$\Ha = \sqrt{R_{\nu}R_mS_c}$ is the Hartmann number, $\f=\g=\boldsymbol{0}$, and boundary conditions are given by
\begin{align*}
& \u = \boldsymbol{0} \quad \mathrm{on}\; x_2= \pm 1 ,\\
& (p\I  - R_{\nu}^{-1}\nabla \u)\n = p\n \quad \mbox{on}\; x_1=0 \; \mbox{and}\; x_1=6,\\
& \n \times \b = \n \times \b_d \quad \mbox{on}\; \partial \O.
\end{align*}

Without loss of generality, we set $G=0.1$. The numerical experiment is carried out on meshes $\T_h^2$ with $h=0.358$. Figure \ref{Fig4} presents the first components of numerical solutions ($L^2$-projection) and analytical solutions with $\Ha = 1$ ($R_{\nu}=1,R_m=0.1,S_c=10$) and $\Ha = 5$ ($R_{\nu}=5,R_m = 1,S_c=5$). The experiment indicates that the numerical solutions almost agree with the analytical solutions.
\begin{figure}[H]
\vspace{-0.2cm}
\centering
\subfigure
{
    \begin{minipage}[b]{.4\linewidth}
        \centering
        \includegraphics[scale=0.4]{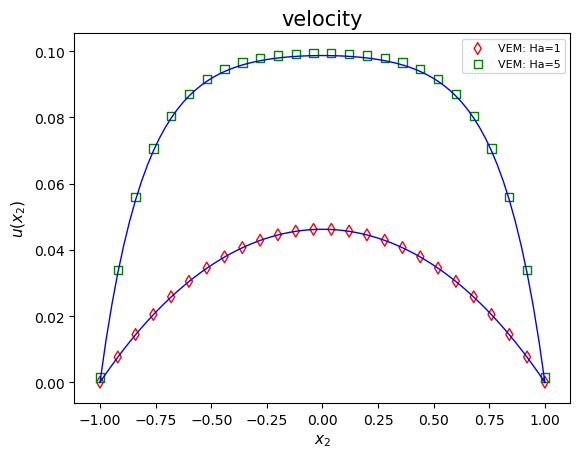}
    \end{minipage}
}
\subfigure
{
 	\begin{minipage}[b]{.4\linewidth}
        \centering
        \includegraphics[scale=0.4]{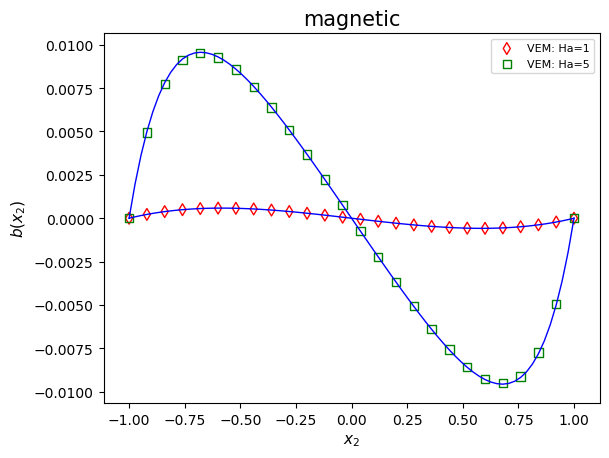}
    \end{minipage}
}
\caption{Numerical (scatter), analytical (line) along $x_1 = 3$.} \label{Fig4}
\end{figure}

\end{example}

\section*{Acknowledgements}
The research was supported by the National Natural Science Foundation of China (Nos: 12071404, 11971410, 12026254), Young Elite Scientist Sponsorship Program by CAST (No: 2020QNRC001), the Science and Technology Innovation Program of Hunan Province (No: 2024RC3158), Key Project of Scientific Research Project of Hunan Provincial Department of Education (No: 22A0136), Postgraduate Scientific Research Innovation Project of Hunan Province (No: CX20230614).

\bibliographystyle{abbrv}
\bibliography{references}

\end{document}